 \documentstyle{article}
 \title{The Deformation of Lagrangian Minimal Surfaces in   Kahler-Einstein
Surfaces}
   \author{Yng-Ing Lee}
\begin{document}
\maketitle

 \newtheorem{lemma}{Lemma} 
\newtheorem{proposition}{Proposition} 
\newtheorem{theorem}{Theorem }
\newtheorem{corollary}{Corollary}
\newtheorem{definition}{Definition}
 
 \noindent A Kahler manifold can  be viewed both as a symplectic manifold 
 and a Riemannian
manifold. These two structures are related by the Kahler form. One can study 
the Lagrangian minimal submanifolds which are Lagrangian with respect to the 
symplectic structure and are minimal with respect to the Riemannian structure.
Lagrangian minimal submanifolds have many nice properties and have been studied
by several authors (see \cite{B}, \cite{CT}, \cite{HL}, \cite{Lee},
\cite{Mc},  \cite{SW}, \cite{SYZ}, \cite{Wo},  \cite{Wo2} etc.).
There are obstructions to the existence of the Lagrangian minimal submanifolds 
in a general Kahler manifold \cite{B}. These obstructions 
do not occur in a Kahler-Einstein manifold. But 
even in this case, the general existence is still unknown.
Most of the discussions of the paper are on compact manifolds without
boundary. We assume this from now on unless other conditions are
indicated. The main result of this paper is the following:
     \vskip.4cm
\noindent
{\bf Theorem 4 }  
{\it Assume that $(N,g_{0})$ is a Kahler-Einstein surface with the 
first Chern class  negative.
Let $[A] $ be a class in the second homology group $  H_{2}(N,Z)$, which
can be represented by a finite union of branched Lagrangian  minimal 
surfaces with 
respect to the metric $g_{0}$. Then with respect to any other metric in 
the connected 
component of
$g_{0}$ in the moduli space of Kahler-Einstein metrics, the class $[A]$ can also be 
represented 
by 
a finite union of branched Lagrangian minimal surfaces.}
   \vskip.4cm
\noindent Note that the complex structure on $N$ is allowed to change 
accordingly.
 An immersed Lagrangian minimal submanifold in a  Kahler manifold with 
negative Ricci curvature is strictly stable (\cite{C}, \cite{MW}, \cite{O}). 
Thus one expects to have
a result as the theorem. However,  
there are some major difficulties  due to the occurrence
of branched points   to realize this expectation.
  In this introduction we   first explain  how the ideas work out in the 
  local deformation of the immersed  case. Then we point out
the difficulties in the  branched case and  how we  solve the problems.
When the Lagrangian minimal surface is immersed,  it is strictly stable 
and thus the Jacobi operator is invertible. By the implicit function theory,
 we can find 
a  minimal surface for any nearby metric and the tangent bundle of the surface changes smoothly. Hence the
minimal surface obtained is totally real if the metric is sufficiently close to the original one.
A totally real (branched) minimal surface in a Kahler-Einstein surface  with 
negative scalar curvature  is  Lagrangian (\cite{CT}, \cite{Wo}).
Therefore, we get the local deformation of an immersed
 Lagrangian minimal surface. If we want to continue the process, we need to
take the limit
of a sequence of surfaces
and it is not enough to restrict to the immersed case. We need to extend each 
step to the branched
case. It seems that there is no known result  for the
deformation of  branched minimal surfaces except the holomorphic curves.
The Jacobi operator on a branched minimal 
surface is degenerate and it is a delicate problem to decide the allowable
variations. For our problem, it is certainly not 
enough to consider  
only the variations with support away
from the branched points. Branched minimal immersions are the critical
points of the energy functional.
     We thus study
the  problem  in the map settings and  show that the strict stablity in the sense of Definition~1 works suitably
 for the deformation of a branched minimal immersion. In particular, we have:
    \vskip.4cm
\noindent
 {\bf Theorem 2 }
{\it Assume that  $ \varphi_{0}:\Sigma \rightarrow (N^{n},g_{0})$ is a
strictly stable branched
 minimal immersion. Then there  exists a strictly stable branched
  minimal 
 immersion 
 $ \varphi_{t}:\Sigma \rightarrow (N^{n},g_{t})$  
  for any $g_{t}$  which is close   to $g_{0}$. Furthermore, $\varphi_t$
converges to 
  $\varphi_0$ in $C^{\infty}$ if $g_t$
converges to 
  $g_{0}$ in $C^{\infty}$. }
    \vskip.4cm
\noindent
Here $\Sigma$ is a closed surface and $N^n$ is a complete Riemannian
$n$-manifold which is not necessarily compact.
  We show that a branched Lagrangian
minimal surface in a Kahler surface with negative Ricci curvature is
strictly stable. Thus we can deform the branched Lagrangian
minimal surface to get a family of strictly stable
branched minimal surfaces. We still need to show that these surfaces are 
Lagrangian.
One can hardly control the behavior of the tangent 
bundles after perturbing the branched points.  The perturbation of the holomorphic
curve  $(z^2,z^3)$ reveals some of the complexity. However, there are still
some control in the holomorphic case. We prove that 
when the branched minimal surfaces are stable, we still have
the similar control. More precisely, we show:
   \vskip.4cm
\noindent
 {\bf Theorem 3 }
{\it Let $\varphi_{i} :\Sigma \rightarrow (N,g_i )$ be  a stable 
branched minimal immersion   from a closed  surface $\Sigma$ to a Riemannian  $4$-manifold $(N,g_i )$. Assume that $g_i$ converges to $g_0$ and $\varphi_{i}$ 
converges to 
   $\varphi_0$ in $C^{\infty}$, where $\varphi_0$ is a branched minimal 
immersion 
   from $\Sigma$ to $(N,g_0 )$.
 Then $$a(\varphi_{0} (\Sigma))=\lim_{i \rightarrow \infty} a(\varphi_{i} 
 (\Sigma)).$$}
\vskip.1cm
\noindent The adjunction number $a(\varphi_{i} (\Sigma))$ in the theorem is defined
to be the sum of the integral of the Gaussian curvature on the tangent bundle
and the integral of the Gaussian curvature on the normal bundle. 
It is equal to the total  
 number of the complex points with indices when $N$ has an almost complex 
structure and the complex points on $\varphi_{i}(\Sigma) $ are isolated
\cite{CT}.
The immersed version of the theorem is proved by J.~Chen and G.~Tian \cite{CT}
using a different approach. From this theorem we can conclude that the  
 branched minimal surface obtained above
is totally real when the metric is sufficiently close to the original one.
Thus it is  Lagrangian (\cite{CT}, \cite{Wo}). This shows the local deformation of a branched 
Lagrangian minimal surface. The rest of the proof for {Theorem~4} follows
from an area bound and standard arguments. 
 \newcommand{\VJ}{\vskip .4cm}
\newcommand{\NO}{\noindent}

\VJ

\noindent The organization of the paper is as follows. In section 1 we  study
the critical points of the energy functional and the stablity. 
 This point of view helps us to understand the branched minimal immersions
and the results in this section should have their own interest. 
The
local deformation of a strictly stable branched minimal immersion  is obtained
in section 2. 
We use the three circle  theorem in
section 3 to study the oscilation of the conformal   harmonic maps.  
    The adjunction
number and some necessary preliminaries  are introduced in section 4. 
In section 5  we prove 
the theorem   about the limit of the adjunction numbers. In the last section 
we complete
the proof of the main theorem and give one application.

\VJ

\NO The author would like to thank the interest of G.~Tian and S.~T.~Yau in this work. 
She also wants to express her gratitude to J.G.~Wolfson's insipring discussions
and interest.   Special thanks are to R.~Schoen who brought this problem to the
author's attention and shared with her many of his profound insights. She
is indebted to the refree for the valuable comments and suggestions.
 During the preparation
of this work, the author were supported by  the National Science Council
of Taiwan under two projects and a visit to Stanford University\footnote{partially supported by NSC-84-2121-M-002-008, NSC-87-2115-M-002-012 and 36128F.}. She would like
to thank these organizations as well.
 
\section{The energy functional}
Let $\Sigma$
be a closed surface of genus $r$ and $(N^{n},g)$ be 
a complete Riemannian
$n$-manifold which is not necessarily compact.
The energy functional   on  $C^{\infty}
 (\Sigma, N)\times{\cal M}(\Sigma) \times{\cal M}(N)$ is defined  to be
  $$E(\varphi,h,g)=\int_{\Sigma}\sum
h^{ij} 
  g_{kl} 
\frac{\partial \varphi  ^{k}}{\partial x^i}
\frac{\partial \varphi  ^{l}}{\partial x^j}\;dA,$$ where $C^{\infty}
 (\Sigma, N)$ is the set of smooth maps between $\Sigma$ and $N$,
 $ {\cal M}(\Sigma)$
 and ${\cal M}(N)$ are the set of smooth metrics on $\Sigma$ and $N$ 
respectively, and   $dA$ is the
volume form of $h$ on $\Sigma$. We will use the convention that two same indices
in a term indicate a summation. The quantity $$ \sum h^{ij}   g_{kl}  
\frac{\partial \varphi  ^{k}}{\partial x^i}
\frac{\partial \varphi  ^{l}}{\partial x^j}$$ is denoted by $e(\varphi)$,
which is called the energy density of $\varphi$ between 
$(\Sigma, h )$ and $(N^{n},g)$. If we fix $h,\;g$ and vary the map 
$\varphi$ only, a critical point of $E( \cdot,h, g)$
 is called
a harmonic map between $(\Sigma, h )$ and $(N^{n},g)$. 
 There has been
a thorough study on harmonic maps. Here we only refere to \cite{EL2}, \cite{ES},
\cite{SY2} and the 
reference therein. If we fix  $\varphi, \;g$ and vary the metric $h$ only,
a direct computation gives the following formula:
\begin{lemma}
Assume that $h_t$ is a smooth family of metrics on $\Sigma$
with $h_{0}=h$. Then
\begin{eqnarray*}
\frac{d E(\varphi,h_{t}, g)}{dt}|_{t=0}
 =\int\sum(\frac{1}{2}h^{ij}h^{ \alpha \beta}-h^{i\alpha }h^{ \beta j})\varphi ^{*}(g)_{ij} 
\dot{h}_{\alpha \beta}\; dA,\end{eqnarray*}
where $$\varphi ^{*}(g)_{ij}=\sum   g_{kl}  
\frac{\partial \varphi  ^{k}}{\partial x^i}
\frac{\partial \varphi  ^{l}}{\partial x^j}$$ is the pull back metric and $\dot{h}=\frac{\partial 
h_t}{\partial t}|_{t=0}.$ In particular, a critical point of $E(\varphi, \cdot, g)$ satisfies $\varphi ^{*}(g) =\frac{1}{2}e(\varphi)h$.
That is, the map 
$\varphi$ is weakly conformal.
\end{lemma}\NO {\bf Proof:}

 \medskip
\NO Assume that $x^{1}, x^{2}$ are the local
coordinates on $\Sigma$. Then $dA_{t}=\sqrt{\det h_{t}}\;dx^{1}dx^{2}$  and the energy can be written as $$E(\varphi,h_{t},g)=\int \sum
h_{t}^{ij}\varphi ^{*}(g)_{ij}\sqrt{\det h_{t}}\;dx^{1}dx^{2}  .$$
Since $$\frac{\partial (h_{t}^{ij} \sqrt{\det h_{t}})}{\partial t}
=-\sum h_{t}^{i\alpha}(\dot{h}_{t})_{ \alpha \beta}h_{t}^{\beta j}\sqrt{\det h_{t}}+\sum h_{t}^{ij}\frac{1}{2}h_{t}^{\alpha \beta } (\dot{h}_{t})_{\alpha \beta}\sqrt{\det h_{t}},$$
 the formula follows. If $\frac{d E(\varphi,h_{t}, g)}{dt}|_{t=0}=0$ for
   arbitrary $\dot{h}_{\alpha \beta}$, one has
 $$\sum_{i,j}(\frac{1}{2}h^{ij}h^{\alpha \beta}-h^{i\alpha }h^{ \beta j})\varphi ^{*}(g)_{ij}=0,
\hskip.3cm \mbox{for any}   \hskip.3cm \alpha ,\beta.$$
In the matrix  expression, this becomes $$\frac{1}{2}e(\varphi)h^{-1}-h^{-1}\varphi ^{*}(g)h^{-1}=0.$$ Hence  $\varphi ^{*}(g) =\frac{1}{2}e(\varphi)h$.
\begin{flushright}{\bf Q.E.D.}\end{flushright}

\begin{lemma}
Follow the  notation as in   
Lemma~1 and assume that $h$ is a critical point of $E(\varphi,\cdot, g)$. 
Then  $$\frac{d^{2} E(\varphi,h_{t}, g)}{dt^{2}}|_{t=0}
=\int(\frac{1}{2}e(\varphi)(T\!r\, h^{-1}\dot{h}h^{-1}\dot{h})-
\frac{1}{4}e(\varphi)(T\!r\,h^{-1}\dot{h})^{2})\;dA,$$ where $T\!r$ denotes
the trace of a matrix.
\end{lemma}
\NO {\bf Proof:}

 \medskip
\NO Now we continue the computation in the proof of Lemma~1 and differentiate
$\frac{d E(\varphi,h_{t}, g)}{dt}$. Because
$$\sum (\frac{1}{2}h^{ij}h^{\alpha \beta}-h^{i\alpha }h^{ \beta j})\varphi ^{*}(g)_{ij}=0,$$ 
 those terms which come from the differentiation on $(\dot{h}_{t})_{\alpha \beta}\sqrt{\det h_{t}}$  has no contribuition. So we only need to 
differentiate 
$$\sum (\frac{1}{2}h_{t}^{ij}h_{t}^{\alpha \beta}-h_{t}^{i\alpha }h_{t}^{ \beta j})\varphi ^{*}(g)_{ij} . $$ Now we compute the formula in terms of matrices and get 
\begin{eqnarray*}
&&\frac{d}{dt}(\frac{1}{2}(T\!r\,h^{-1}_{t}\varphi ^{*}(g))h^{-1}_{t}-
h^{-1}_{t}\varphi ^{*}(g)h^{-1}_{t})|_{t=0}\\&=&-\frac{1}{2}(T\!r\,h^{-1}\dot{h}h^{-1}
\varphi ^{*}(g))h^{-1}- \frac{1}{2} e(\varphi) h^{-1}\dot{h}h^{-1} + h^{-1}\dot{h}h^{-1}\varphi ^{*}(g)h^{-1} +h^{-1}\varphi ^{*}(g)h^{-1}\dot{h}h^{-1}\\&=&-\frac{1}{4}e(\varphi)(T\!r\,h^{-1}\dot{h})h^{-1}
+\frac{1}{2} e(\varphi) h^{-1}\dot{h}h^{-1}.
\end{eqnarray*}
Therefore, 
 $$\frac{d^{2} E(\varphi,h_{t}, g)}{dt^{2}}|_{t=0}
=\int(\frac{1}{2}e(\varphi)(T\!r\, h^{-1}\dot{h}h^{-1}\dot{h})-
\frac{1}{4}e(\varphi)(T\!r\,h^{-1}\dot{h})^{2})\;dA.$$
\begin{flushright}{\bf Q.E.D.}\end{flushright}

\begin{lemma}
Assume that  $(\varphi, h)$ is  a critical point of $ E(\cdot,\cdot, g)$.
Let  $h_t$ be a smooth family of metrics on $\Sigma$
with $h_{0}=h,\; \dot{h}=\frac{\partial 
h_t}{\partial t}|_{t=0}$ and let $  \varphi_{t} $
be a smooth family of maps from $\Sigma$ to $N$ with $\,\varphi_{0} =
\varphi, \; \frac{\partial \varphi_{t}}{\partial t}|_{t=0}=V$. Then we have
\begin{eqnarray*}
&&\frac{d^{2} E(\varphi_{t},h_{t}, g)}{dt^{2}}|_{t=0}\\&=&
2\int (|\nabla V|^{2}+\sum
<R^{N}(d\varphi(e_{i}),V)d\varphi(e_{i}),V>)\;dA\\&&
+\int(\frac{1}{2}e(\varphi)(T\!r\, h^{-1}\dot{h}h^{-1}\dot{h})-
\frac{1}{4}e(\varphi)(T\!r\,h^{-1}\dot{h})^{2})\;dA\\&&
+2\int\sum(\frac{1}{2}h^{ij}h^{\alpha \beta}-h^{i\alpha}h^{\beta j})  
\dot{h}_{\alpha \beta}(<\nabla _{\frac{\partial}{\partial x^{i}}}V,d\varphi(\frac{\partial}{\partial x^{j}}) >+
<\nabla _{\frac{\partial}{\partial x^{j}}}V,d\varphi(\frac{\partial}{\partial x^{i}}) >)\;dA,
\end{eqnarray*}
where $R^{N}$ is the curvature tensor of $(N,g)$ and $\{e_{1},\,e_{2}\}$ is a local frame for $h$.
\end{lemma}
\NO {\bf Proof:}

 \medskip
\NO By Leibniz's rule, one has 
\begin{eqnarray*}
\frac{d^{2} E(\varphi_{t},h_{t}, g)}{dt^{2}}|_{t=0}= \frac{d^{2} E(\varphi_{t},h, g)}{dt^{2}}|_{t=0}+\frac{d^{2} E(\varphi,h_{t}, g)}{dt^{2}}|_{t=0}+2\frac{\partial ^{2} E(\varphi_{t},h_{s}, g)}
{\partial t \partial s}|_{t=s=0}
\end{eqnarray*}
The formula $$\frac{d^{2} E(\varphi_{t},h, g)}{dt^{2}}|_{t=0}=2\int (|\nabla V|^{2}+\sum
<R^{N}(d\varphi(e_{i}),V)d\varphi(e_{i}),V>)\;dA$$ is well known
and can be found for instance in \cite{EL2} or \cite{SY2}. The formula for
$\frac{d^{2} E(\varphi,h_{t}, g)}{dt^{2}}|_{t=0}$ is derived in Lemma~2.
A direct computation shows that $$\frac{\partial \varphi_{t} ^{*}(g)_{ij}}{\partial t}=\, <\nabla _{\frac{\partial}{\partial x^{i}}}V,d\varphi(\frac{\partial}{\partial x^{j}}) >+
<\nabla _{\frac{\partial}{\partial x^{j}}}V,d\varphi(\frac{\partial}{\partial x^{i}}) >.$$  This together with the computation in the proof of Lemma~1 give
the formula for $\frac{\partial ^{2} E(\varphi_{t},h_{s}, g)}
{\partial t \partial s}|_{t=s=0}$.
\begin{flushright}{\bf Q.E.D.}\end{flushright}
 The energy is a conformal  
 invariant on the metric of the domain when the domain is two dimensional. Denote $E_{g}(\varphi,h)=E(\varphi,h,g).$ Then $E_{g}$ can be viewed as a 
smooth function $E_{g}(\varphi, [h])$
 on $C^{\infty} (\Sigma, N)\times{\cal T}_r$, where ${\cal T}_r$ is the 
Teichmuller
 space of $\Sigma$ and $[h]$ is the conformal class of $h$. The branched immersions were defined and studied 
in \cite{GOR}. In particular, the pull back metric  $\varphi ^{*}(g)$ of 
a branched immersion
can be expressed as 
$\lambda 
^{2}
 h,$ where $h$ is a smooth metric on $\Sigma$ and $\lambda$ is a smooth scalar function with  isolated and
finite order zeros. The zeros of $\lambda$ are called  the branched points
of $\varphi$. Hence in this case $[\varphi ^{*}(g)]$
is well-defined   and $E_{g}(\varphi, [\varphi ^{\ast} (g)])
=2 A(\varphi, g)$, where $A(\varphi, g)$ is the area of $\varphi(\Sigma)$ in 
$(N,g)$. 
\VJ
\NO{\bf Remark:} Assume that  $(\varphi, [h])$ is  a critical point of $E_{g}$ and $[h_{t}]$
is a variation of the conformal structure.    
Because in our case the energy functional is a conformal invariant and by a result of 
Moser \cite{Mo}, we can choose $h_{t}$ such that they all determine the
same volume form. 
That is, we can assume $T\!r\,h^{-1}\dot{h}=0$ in the second variational formula. 
\VJ
\NO One thus has 
  the following result of J.~Sacks and K.~Uhlenbeck:
\begin{corollary}
\cite{SU} The map $ \varphi$ is a branched minimal immersion if and only if
 $(\varphi, [h])$ is  a critical point of $ E_{g}(\cdot,\cdot)$ and  
 $\varphi $ is 
a nonconstant map,  where the smooth metric $h$ is conformal to the pull back metric $\varphi ^{*}(g)$.
\end{corollary}
\NO {\bf Proof:}

 \medskip

\NO If $(\varphi, [h])$ is  a critical point of $ E_{g}(\cdot,\cdot)$, certainly
$\varphi $ is  a critical point of $ E(\cdot, h , g)$ and $h$ is
 a critical point of $ E(\varphi,\cdot , g)$. Thus the map $\varphi$ is both harmonic and weakly conformal between $(\Sigma, h )$ and $(N^{n},g)$.
When $\varphi $ is a nonconstant map, this exactly means that $\varphi$ is a branched minimal immersion. 
If $ \varphi$ is a branched minimal immersion, then $\varphi$ is both harmonic and weakly conformal between $(\Sigma, h )$ and $(N^{n},g)$. Thus $\varphi $ is  a critical point of $ E(\cdot, h , g)$ and $h$ is
 a critical point of $ E(\varphi,\cdot , g)$.
Assume that  $h_t$ is a smooth family of metrics on $\Sigma$
with $h_{0}=h,\; \dot{h}=\frac{\partial 
h_t}{\partial t}|_{t=0}$ and   $  \varphi_{t} $
is a smooth family of maps with $\,\varphi_{0} =
\varphi, \; \frac{\partial \varphi_{t}}{\partial t}|_{t=0}=V$.
By Leibniz's rule, one has 
\begin{eqnarray*}
\frac{d E(\varphi_{t},h_{t}, g)}{dt}|_{t=0}= \frac{d E(\varphi_{t},h, g)}{dt}|_{t=0}+\frac{d E(\varphi,h_{t}, g)}{dt}|_{t=0}=0.
\end{eqnarray*}
Hence $(\varphi, [h])$ is  a critical point of $E_{g}(\cdot,\cdot)$.  \begin{flushright}{\bf Q.E.D.}\end{flushright}

\VJ
\NO Assume that $ \varphi
:\Sigma \rightarrow (N^{n},g)$ is a branched minimal immersion
and $(\varphi, [h])$ is  the corresponding critical point on $E_{g}$.  
Define a function $f_\varepsilon $ on $\Sigma$:
\begin{eqnarray}
f_\varepsilon (x)=\left \{ \begin{array}{cll}
0&\hspace{2.5em}\mbox{\ }& |x|< \varepsilon^{2}\\ 
\displaystyle{{\log \frac{|x|}{\varepsilon^{2}}}\over
{\log \frac{1}{\varepsilon }}}
&\mbox{     } &\varepsilon^{2}\leq |x|\leq \varepsilon\\ 
1&
\mbox{     }& |x|> \varepsilon .
\end{array} \right.
\end{eqnarray}
  Then $\lim\limits_{\varepsilon  \rightarrow 0}\int  |\nabla f_\varepsilon|^{2} \,d A=0.$  Now we choose $f_\varepsilon $ such that it vanishes near each branched point of $\varphi$.
\begin{lemma}
If we denote the second variation of $E_{g}$ in the direction of $V$ and $\dot{h}$ by $\delta^{2}E_{g}(V, \dot{h})$, then
$$\delta^{2}E_{g}(V, \dot{h})= \lim_{\varepsilon  
\rightarrow 0}
 \delta ^{2} E_{g}(f_\varepsilon V,\dot{h})$$
\end{lemma}
\NO {\bf Proof:}

 \medskip
\NO A direct computation gives us $$\nabla_{ \frac{\partial}{\partial x^i}}
   f_\varepsilon  V
  =f_\varepsilon \nabla_{ \frac{\partial}{\partial x^i}}  V
  +\frac{\partial f_\varepsilon }{\partial x^i}V $$
 and $$|\nabla  f_\varepsilon  V|^{2}= f_\varepsilon ^{2}
  |\nabla  V|^{2} +|\nabla f_\varepsilon|^{2}|V |^{2}
  +2\sum  <e_{i}( f_\varepsilon )V,
  f_\varepsilon \nabla_{ e_{i}}  V>.$$
Since $h^{-1},\, \dot{h},\, \varphi,$ and $V$ are smooth and fixed, the norms and the
norms of their derivatives are all bounded. Therefore,
\begin{eqnarray*}
&&|\delta^{2}E_{g}(V, \dot{h})- \lim_{\varepsilon  
\rightarrow 0}
 \delta ^{2} E_{g}(f_\varepsilon V,\dot{h})|\\&& \leq C_{1}
\lim_{\varepsilon  \rightarrow 0}\int  |\nabla f_\varepsilon|^{2} \,d A
+C_{2}
(\lim_{\varepsilon  \rightarrow 0}\int  |\nabla f_\varepsilon|^{2} \,d A)^{\frac{1}{2}}\\&&=0,  \end{eqnarray*} 
 where $C_{1}$ and $C_{2}$ are positive constants
independent of $\varepsilon$.
\begin{flushright}{\bf Q.E.D.}\end{flushright}
\begin{definition}
A   branched minimal immersion
$ \varphi  :\Sigma \rightarrow (N^{n},g)$ is called
strictly stable    if $\lim\limits_{\varepsilon  \rightarrow 0}\delta ^{2} A(f_{\varepsilon}V)>0$ for any  $V=\frac{\partial \varphi _{t}}{\partial t}|_{t=0}$, where $f_{\varepsilon}$ is chosen as in (1) and $  \varphi_{t} $
 is   a smooth family of maps from $\Sigma$ to $N$   with $\varphi_{0} =
\varphi$. It is called
stable if $\lim\limits_{\varepsilon  \rightarrow 0}\delta ^{2} A(f_{\varepsilon}V)\geq 0$
\end{definition} 
 
\begin{theorem}
A branched minimal immersion $ \varphi  :\Sigma \rightarrow (N^{n},g)$ is strictly stable if and only if
the correponding critical point on $E_{g}$ is strictly stable.
\end{theorem}

\NO {\bf Proof:}

 \medskip
\NO We first claim that for any branched immersion
 $ \phi  :\Sigma \rightarrow (N^{n},g)$  and any smooth metric $h$ on $\Sigma$,
one always has  $$E_{g}(\phi,h)\geq   2 A(\phi, g).$$
Choose $x^{1}, x^{2}$ to be the conformal coordinates  for the pull back
metric $\phi ^{\ast} (g)$. That is, $$\sum   g_{kl} 
\frac{\partial \phi  ^{k}}{\partial x^i}
\frac{\partial \phi  ^{l}}{\partial
x^j}=\mu ^{2}\delta _{ij},$$ where $\mu$ is nonnegative.
Express the inverse matrix $(h^{ij})$ in this coordinates as 
$(\begin{array}{ll} a&c\\c&b \end{array}) $, where $a$ and $b$ are
positive. Then we have


\begin{eqnarray*}
E_{g}(\phi,h)&=& \int\sum h^{ij} 
  g_{kl} 
\frac{\partial \phi  ^{k}}{\partial x^i}
\frac{\partial \phi  ^{l}}{\partial
x^j}\;dA\\ 
&=&\int(a\mu+b\mu)
\frac{1}{\sqrt{ab-c^2}}dx^{1}dx^{2}\\ 
&\geq& 2\int \mu 
\frac{\sqrt{ab } }{\sqrt{ab-c^2}}dx^{1}dx^{2}\\ 
&\geq&  2 A(\phi, g).\end{eqnarray*} 
The equalities hold if and only if $a=b$ and $c=0$, i.e.,
when $\phi$ is a weakly conformal   map.
\VJ
\NO Assume that   $(\varphi, [h])$ is  the corresponding critical point on $E_{g}$  of the branched minimal immersion,
 where $h$ is a smooth metric  on $\Sigma$. Let  $h_t$ be a smooth family of metrics on $\Sigma$
with $h_{0}=h,\; \dot{h}=\frac{\partial 
h_t}{\partial t}|_{t=0}$ and $  \varphi_{t} $
be a smooth family of maps from $\Sigma$ to $N$  with $\,\varphi_{0} =
\varphi, \; \frac{\partial \varphi_{t}}{\partial t}|_{t=0}=V$.  
Define $\varphi_{t}^{\varepsilon}(x)=\exp_{\varphi (x)} tf_{\varepsilon} V(x)$,
where $f_{\varepsilon}$ is chosen as in (1).
 Then  $\varphi_{t}^{\varepsilon}$ is  a smooth family of branched immersions
from $\Sigma$ to $N$ with $\,\varphi_{0}^{\varepsilon} =
\varphi$ and $ \frac{\partial \varphi^{\varepsilon}_{t}}{\partial t}|_{t=0}
=f_{\varepsilon}V$. 
By the claim proved above, one has $$E_{g}(\varphi_{t}^{\varepsilon}, h_{t})
\geq 2 A(\varphi_{t}^{\varepsilon},g).$$  Define the $C^{2}$ nonnegative
  function $F$ by $$F(t)=E_{g}(\varphi_{t}^{\varepsilon}, h_{t})
- 2 A(\varphi_{t}^{\varepsilon},g).$$ Because $F(0)=\dot{F}(0)=0$, it follows that $\ddot{F}(0)
\geq 0.$ Hence $$\delta ^{2} E_{g}(f_\varepsilon V,\dot{h}) \geq 2 \delta^{2}
A(f_\varepsilon V)$$  and thus $$\delta ^{2} E_{g}(  V,\dot{h}) \geq \lim_{\varepsilon  \rightarrow 0} 2 \delta^{2}
A(f_\varepsilon V)>0.$$ One also has $\delta ^{2} E_{g}(0,\dot{h})>0$ by Lemma~3 for $\dot{h}$ which is not identically zero. This shows that  $(\varphi, [h])$ is  a strictly stable critical point on $E_{g}$.
\VJ
\NO Assume that $(\varphi, [h])$ is  a strictly stable critical point on $E_{g}$. Then $ \varphi$ is a branched minimal immersion and one has $$
\delta^{2}
A(f_\varepsilon V)=\frac{1}{2}\delta ^{2} E_{g}(f_\varepsilon V,0).$$ Thus 
$$\lim_{\varepsilon  \rightarrow 0} \delta^{2}
A(f_\varepsilon V)=\frac{1}{2}\delta ^{2} E_{g}( V,0)>0.$$ Hence
$ \varphi$ is a strictly stable branched minimal immersion
\begin{flushright}{\bf Q.E.D.}\end{flushright}
\section{The deformation of branched minimal surfaces}
Let $\Sigma$
be a closed surface and $(N^{n},g_{0})$ be a complete Riemannian
$n$-manifold which is not necessarily compact. The strict stablity 
in the sense of Definition~1 works suitably
 for the deformation of a branched minimal immersion. In particular, we have:
\begin{theorem}
 Assume that  $ \varphi_{0}:\Sigma \rightarrow (N^{n},g_{0})$ is a
strictly stable branched
 minimal immersion. Then there  exists a strictly stable branched
  minimal 
 immersion  $ \varphi_{t}:\Sigma \rightarrow (N^{n},g_{t})$ for any $g_{t}$ which is close enough to $g_{0}$. Furthermore, $\varphi_t$
converges to 
  $\varphi_0$ in $C^{\infty}$ if $g_t$
converges to 
  $g_{0}$ in $C^{\infty}$.
 \end{theorem}


\NO {\bf Proof:}

 \medskip

 \NO Let  $(\varphi_{0}, [h_{0}])$ be  the corresponding critical point on $E_{g_{0}}$. By Theorem~1, one knows that $(\varphi_{0}, [h_{0}])$
is  strictly stable. Particularlly,  $\varphi_{0}$
is a strictly stable harmonic map from $(\Sigma,  {h}_{0})$ to $(N,g_{0})$.
It is a theorem of Eells and Lemaire   \cite{EL} that there exists a
neighborhood $\cal V$ of $ h_{0}$ and $g_{0}$ in ${\cal M}(\Sigma)
\times {\cal M}(N)$  and a unique smooth map $S$ on $\cal V$ such that 
$S(  {h}_{0}, g_{0})=\varphi _{0}$ and $S(h,g)$ is a smooth harmonic map 
between $(\Sigma, h)$ and $(N, g)$. Let $\varphi _{t,h}=S(h,g_{t})$ and 
$\cal U$ be the corresponding neighborhood of $[h_{0}] $ in
the 
Teichmuller
 space ${\cal T}_r$. Since the energy is a conformal invariant on the domain,
$\varphi _{t,h}$ is also harmonic with respect to any other representative
of $[h]$. Thus $\varphi _{t,h}$ is determined by $[h]$ in ${\cal U}.$
Define $\bar{E}_{g_t}:{\cal U} \rightarrow R$ by $\bar{E}_{g_t}([h])=E_{g_t}
(\varphi _{t,h},[h]).$ Then $[h]$ is a critical point of $\bar{E}_{g_t}$ if and only if $(\varphi _{t,h},[h])$ is a critical point of $E_{g_t}(\cdot,\cdot).$
The differential $d\bar{E}_{g_t}|_{[h]}$ lies in $T^{\ast}_{[h]}{\cal T}_r$. It is identified with $R^{6r-6}$ if we choose local coordinates near $[h_0]$.
Define $$G: {\cal U} \times (-\varepsilon, \varepsilon)
 \rightarrow d\bar{E}_{g_t}|_{ [h]}.$$
  We have that 
 $G([h_{0}],0)=0$ and $dG|_{([h_{0}],0)}$ is of full rank because     $(\varphi_{0},[h_{0}])$ is a strictly stable critical point of $E_{g_0}$.
  By applying the implicit function theory to $G$, there exists a smooth
  path $[h_{t}]$ in ${\cal T}_{r}$ such that $G([h_{t}],t)=0$. Hence $[h_t]$
is a critical point of $\bar{E}_{g_t}$. Denote 
$\varphi _{t,h_{t}}$ by $\varphi _{t}$. It follows that 
$(\varphi _{t },[h_{t}])$ is a critical point of $E_{g_t}$ and $\varphi _{t}$
is a branched minimal immersion. Because
the energy $E_g$ depends smoothly on $g$, we can conclude that 
$(\varphi _{t },[h_{t}])$ is a strictly stable critical point for $t$ small enough.
  Thus $\varphi_{t}$ is a strictly stable branched
 minimal immersion.  By the construction of $\varphi _{t}$
 and the theorem of Eells and Lemaire in \cite{EL},
one also has 
 $\varphi_t$ converges to 
  $\varphi_0$ in $C^{\infty}$.
\begin{flushright}{\bf Q.E.D.}\end{flushright}
\begin{proposition}
 Every    
branched Lagrangian  minimal immersion in a  Kahler  surface $N$ with negative
Ricci curvature is strictly
stable.
\end{proposition}
{\bf Proof:}

\medskip

\NO    Let $f_{\varepsilon }$ be   defined as in (1), which has support 
away from the branched points of $\varphi _{0}$ and assume that $V$ is a
vector field  along $\varphi _{0}$ which is defined on $\Sigma$.
Define the one 
form $\beta_\varepsilon $ on $\Sigma$ by $\beta_\varepsilon (u)=
<Jf_{\varepsilon }V ,u>$, where $J$ is the complex structure on $N$ and
$u\in T\Sigma$. By a result in
 \cite{C} and \cite{MW}, we have
\begin{eqnarray*}
 \delta ^{2} A(f_\varepsilon V ) 
&=& \int 
_{\Sigma}(|d\beta_{\varepsilon}|^{2}+|\delta 
\beta_{\varepsilon}|^{2}-  Ric\,(f_{\varepsilon }V ,
f_{\varepsilon }V ))\;d A\\  
&\geq& c \int_{\Sigma} |f_{\varepsilon }V |^{2}\;d A,
 \end{eqnarray*}
where $  Ric$ 
 is the Ricci curvature of the Kahler surface  satisfying
$$  Ric\,(V,V)\leq -c|V|^{2}$$ for some positive constant $c$. 
Thus
 $$ \lim_{\varepsilon  \rightarrow 0}
 \delta ^{2} A(f_\varepsilon V)
  \geq \lim_{\varepsilon  \rightarrow 0}
 c\int_{\Sigma} |f_{\varepsilon }V |^{2}\;d A
 >0,$$  and the 
map is strictly stable in  the  sense of Definition~1.
 \begin{flushright}{\bf Q.E.D.}\end{flushright}
 \begin{corollary}
Let $\varphi$ be a branched Lagrangian  minimal immersion in a   Kahler  surface with negative
Ricci curvature. Then there is a strictly stable branched minimal immersion  
 near $\varphi$ with respect to the Riemannian metric which is close  to this 
Kahler metric.
\end{corollary}

\section{The oscillation of the conformal  harmonic maps}  
Let $\varphi:\Sigma \rightarrow
 N$ be a smooth map from a Riemannian surface $\Sigma$ to a complete $n$-dimensional Riemannian manifold
 $N$. Let $\theta ^{1},\; \theta ^{2}$ be an orthonormal coframe in a 
neighborhood of $p \in \Sigma$ and let $\omega^{1},\cdots,\omega^{n}$
be  an orthonormal coframe in a 
neighborhood of $\varphi(p) \in N$. Define $\varphi _{\alpha}^{l},\;
1\leq \alpha \leq 2,\; 1\leq l \leq n $ by
$$\varphi ^{\ast}\omega^{l}=\sum \varphi^{l}_{\alpha}
\theta ^{\alpha} \hskip.3cm \mbox{for} \hskip.3cm 1\leq l \leq n.$$
We have the structure equations for $N$ and $\Sigma$
$$d\omega^{l}=\sum \omega_{m}^{l} \wedge \omega^{m} \hskip.3cm
 \mbox{and}
\hskip.3cm \omega_{m}^{l}=-\omega _{l}^{m} \hskip.3cm \mbox{for} \hskip.3cm
 1\leq l,m \leq n,$$
$$d\theta^{\alpha}=\sum \theta^{\alpha}_{ \beta} 
\wedge \theta^{\beta} \hskip.3cm \mbox{and}
\hskip.3cm \theta^{\alpha}_{ \beta}=-\theta^{\beta}_{ \alpha }\hskip.3cm
\mbox{for} \hskip.3cm 1\leq \alpha ,\beta\leq 2.$$
Define $\varphi^{l}_{\alpha \beta},\;
1\leq \alpha, \beta \leq 2,\; 1\leq l \leq n $ by
$$d \varphi ^{l}_{\alpha}+\sum \varphi^{m}_{\alpha}\varphi ^{\ast}\omega_{m}^{l}+\sum \varphi^{l}_{\beta}\theta^{\beta}_{
\alpha}
=\sum \varphi^{l}_{\alpha \beta}\theta^{\beta}.$$

\VJ

\NO Choose the local coordinates at $p$ to be $0$ and
let $\rho^{2} (y)$   be the square of the distance   between
$y  $ and  $\varphi (0)$ in $(N,g)$. Then $\rho^{2}(\varphi (x))$ is a function
on $\Sigma$ and 
   $$ \Delta \rho ^{2}(\varphi (x))=2\sum_{\alpha}(\sum_{l} \rho_{l}\varphi^{l}_{\alpha})^{2}+2
 \rho\sum \rho_{kl}\varphi^{k}_{\alpha}
 \varphi^{l}_{\alpha}+2
 \rho\sum \rho_{ l}\varphi^{l}_{\alpha \alpha},$$ 
 where $1 \leq \alpha \leq 2 $ and $1 \leq k,l \leq n.$
   The condition of $\varphi$ to be harmonic is equivalent to $\sum 
  \varphi^{l}_{\alpha \alpha}=0$
 for all $l$. If we choose the normal coordinates $\{y^{1},\ldots,y^{n}\}$ at 
 $\varphi(0)$, we have $$\rho_{l}(y)=
 \frac{y^l}{\rho}\;\; \mbox{  and  } \;\;\rho_{kl}(y)=
 \frac{\delta_{kl}}{\rho}-\frac{y^{k}y^{l}}{\rho^3}-\sum \Gamma^{m}_{kl}
 \frac{y^m}{\rho} .$$
 When $\varphi$ is harmonic, one has $$ \Delta \rho ^{2}=2|\nabla
\varphi|^{2}-2\sum
 \varphi^{m}\Gamma^{m}_{kl}\varphi^{k}_{\alpha}
 \varphi^{l}_{\alpha}. $$ Hence $\rho ^{2}(\varphi (x))$ is a subharmonic function 
on $\Sigma$ when the metric on $N$ is flat.
 A general Riemannian metric satisfies $\Gamma^{m}_{kl}(y)=O(|y|). $  By taking $y=\varphi(x)$, 
it follows that   $\rho ^2 (\varphi (x))$ is  subharmonic
 when $|x|$ is small. 
 Further computation shows that \begin{eqnarray*}
 \Delta \log \rho ^{2}&=&\frac{\Delta \rho ^{2}}{\rho ^{2}}-\frac{|\nabla \rho ^{2}|^2}{\rho ^{4}}
 \\ 
&=&\frac{2|\nabla \varphi|^{2}-2\sum
 \varphi^{m}\Gamma^{m}_{kl}\varphi^{k}_{\alpha}
 \varphi^{l}_{\alpha}}{\rho ^{2}}-\frac{4 \sum_{\alpha}(\sum_{l} \rho_{l}\varphi^{l}_
 {\alpha})^{2}}{\rho ^{2}}.\end{eqnarray*} 
 \begin{lemma}
 Assume that  $\varphi:(B_{2}(0), \sum_{i=1}^{2} (dx^{i})^{2}) \rightarrow (N,g)
 $ is a 
 conformal   harmonic map from a ball of radius $2$ into a normal
neighborhood of $\varphi (0)$ in $N$. Then we have $$\max_{B_{r_2}(0)}\rho ^2 (\varphi (x))
  \leq (\frac{r_2}{r_{1}})^{C} \max_{B_{r_1}(0)}
  \rho ^2 (\varphi (x)) $$ for $0<r_{1}\leq r_{2}\leq  \varepsilon
 \leq 1$, where
 $\varepsilon$ is a constant depending only on the metric $g$ and $C$ is a constant 
  independent of $r_1$ and $r_2$.
 \end{lemma}
{\bf Remark: }The radius 2 in the Lemma does not matter and the main point 
is to have a ball of fixed radius which maps into a normal neighborhood of 
$\varphi (0)$. The constant $\varepsilon$ is chosen such that it is less
than the fixed radius and $\rho ^{2}(\varphi (x))$ is subharmonic on 
$B_{\varepsilon}(0)$.
\vskip.4cm

\noindent {\bf Proof:}

 \medskip

\NO When $\varphi$ is a constant map, the lemma holds trivally. So we assume
that $\varphi$ is a nonconstant map.
Because $\varphi$ is a conformal  map, we have 
 $$\sum  (\varphi^{l}_{1})^{2}=\sum  
 (\varphi^{l}_{2})^{2}=\mu^{2}\;\;\mbox{  and  }\;\;\sum  
   \varphi^{l}_{1} \varphi^{l}_{2}=0 ,$$
where $\mu$ is a smooth and nonnegative scalar function with isolated and finite order zeros.
Hence
  $$d \varphi(\frac{\partial}{\partial x^1})
=\mu e_{1}\;\;\mbox{ and }\;\;d \varphi(\frac{\partial}{\partial x^2})
=\mu e_{2},$$  where $ e_{1}$ and $e_2$
are  orthonormal.
Therefore, 
\begin{eqnarray*}
\Delta \log \rho ^{2}&=&\frac{2|\nabla \varphi|^{2}-2\sum
 \varphi^{m}\Gamma^{m}_{kl}\varphi^{k}_{\alpha}
 \varphi^{l}_{\alpha}}{\rho ^{2}}-\frac{4 \sum_{\alpha}(\sum_{l} \rho_{l}\varphi^{l}_
 {\alpha})^{2}}{\rho ^{2}}\\ 
&=&\frac{4\mu^{2}-4\mu^{2}
\sum ( \nabla \rho \cdot e_{\alpha} )^{2}}{\rho ^{2}}
- \frac{ 2\sum
 \varphi^{m}\Gamma^{m}_{kl}\varphi^{k}_{\alpha}
 \varphi^{l}_{\alpha}}{\rho ^{2}} \\ 
& \geq&  \frac{4 \mu ^{2}-4 \mu ^{2} }{\rho ^{2}}- \frac{ 2\sum
 \varphi^{m}\Gamma^{m}_{kl}\varphi^{k}_{\alpha}
 \varphi^{l}_{\alpha}}{\rho ^{2}} \\ 
& \geq& -\bar c, \end{eqnarray*}
 where the positive constant $\bar c$ depends only on the upper bound of
$|\nabla \varphi|^{2} $
 and  the metric $g$.
 We use the fact that $|\nabla
 \rho|=1$ in the first inequality. 
   A direct computation shows that  $\Delta r^{2}=2$ and $\Delta \log r=0$,
    where $r$ is the distance function on the domain. Define $$F(x)=e^{\frac{{\bar c}}{2}\, r(x)^{2}}
    \rho^{2} (\varphi (x)).$$ Then 
\begin{eqnarray*}\Delta \log F(x)&=&\Delta \log \rho ^{2}+
    \Delta  \frac{{\bar c}}{2}\, r^{2}\\
&\geq& -{\bar c}+{\bar c}\\
&=&0,\end{eqnarray*} 
so that $\log F(x)$ is a subharmonic function.
Define $$M(r)=\max_{\partial B_{r}(0)}
 F(x) =\max_{ B_{r}(0)}
  F(x).$$ Then the function $$ \log F(x)-
 \frac{\log r-\log r_{1}}
 {\log r_{2}-\log r_{1}} \log M(r_{2})-\frac{\log r_{2}-\log r }
 {\log r_{2}-\log r_{1} } \log M(r_{1})$$ is a subharmonic function and has 
 nonpositive values on the circles
 of radius $r_1$ and $r_2$. By applying the maximun principle to the 
 annulus between radius
 $r_{1}$ and $r_2$, we conclude that
 $$\log M(r)\leq \frac{\log r-\log r_{1}}
 {\log r_{2}-\log r_{1}} \log M(r_{2})+\frac{\log r_{2}-\log r }
 {\log r_{2}-\log r_{1} } \log M(r_{1})$$ for $r_{1}<r< r_{2}.$ 
 This means that $\log M(r)$ is a 
 convex function in terms of $\log r$. Since the choice of $r_1$
 and $r_2$ is arbitrary, the conclusion holds for all $0<r<2$. Now we want
 to compute the derivative of $\log M(r)$ with respect to $\log r$ at $r=1$  
and bound it by a constant
 $C$. We have
 $$\frac{d\log M(r)}{d \log r}=\frac{d\log M(r)}{d  r}
 \frac{d r}{d \log r}=\frac{M'(r)}{M(r)}r,$$
where $$M(r)=\max_{\partial B_{r}(0)}
 F(x) =e^{\frac{{\bar c}}{2}\, r^{2}} \max_{\partial B_{r}(0)}\rho^{2} (\varphi (x))$$
and $$M'(r)\leq \max_{\partial B_{r}(0)} |\nabla F(x)|.$$
A direct computation shows that  
\begin{eqnarray*}|\nabla F(x)|&\leq &{\bar c} re^{\frac{{\bar c}}{2}\, r^{2}}\rho^{2} (\varphi (x))
 +2e^{\frac{{\bar c}}{2}\, r^{2} }\rho  (\varphi (x))|\nabla \rho  (\varphi (x))|
\\&\leq& {\bar c}r e^{\frac{{\bar c}}{2}\, r^{2} }\rho^{2} (\varphi (x))
 +2e^{\frac{{\bar c}}{2}\, r^{2} }\rho  (\varphi (x))|\nabla  \varphi (x)| \end{eqnarray*} 
 for $x \in \partial B_{r}(0)$. Hence
\begin{eqnarray*}\frac{M'(1)}{M(1)} \leq {\bar c}+2\frac{\max _{\partial
B_1} 
\rho(\varphi (x)) |\nabla \varphi| }{\max _{\partial
B_1} 
  \rho^{2}(\varphi (x)) }  \leq {\bar c}+2\frac{\max _{\partial
B_1} 
 |\nabla \varphi(x)| }{\max _{\partial
B_1} 
  \rho(\varphi (x)) }. \end{eqnarray*}So we can choose
  $$C={\bar c}+2\frac{\max _{\partial
B_1} 
 |\nabla \varphi(x)| }{\max _{\partial
B_1} 
  \rho(\varphi (x)) }.$$
Because the slope of a convex 
 function is increasing, we have that 
 $$\frac{\log M(r_{2})-\log M( r_{1})}
 {\log r_{2}-\log r_{1}}\leq C,$$ for $0<r_{1}< r_{2}\leq 1.$ Therefore,
$$\log \frac{M(r_{2})}{M(r_{1})}
 \leq C\log \frac{ r_{2}}{ r_{1} },$$ or $$\frac{M(r_{2})}{M(r_{1})}
 \leq(\frac{ r_{2}}{ r_{1} })^{C}.$$ Thus we have 
 $$\frac{ e^{\frac{{\bar c}}{2}\, r^{2}_{2} }\max_{\partial B_{r_2}}\rho^{2} (\varphi (x))}
 { e^{\frac{{\bar c}}{2}\, r^{2}_{1} }\max_{\partial B_{r_1}}\rho^{2} (\varphi (x)) }
 \leq(\frac{ r_{2}}{ r_{1} })^{C}.$$ Choose   $\varepsilon$ such that $\rho^2 (\varphi (x))$
 is subharmonic when $|x|\leq \varepsilon$. Hence 
 $$\max_{\partial B_{r }}\rho^{2} (\varphi (x)) =
\max_{  B_{r }}\rho^{2} (\varphi (x))  $$ for $r\leq \varepsilon$.
It follows that $$\frac{\max_{  B_{r_2 }}\rho^{2} (\varphi (x)) }
{\max_{  B_{r_1 }}\rho^{2} (\varphi (x)) }\leq
\frac{ e^{\frac{{\bar c}}{2}\, r^{2}_{2} }\max_{\partial B_{r_2}}\rho^{2} (\varphi (x)) }
 { e^{\frac{{\bar c}}{2}\, r^{2}_{1} }\max_{\partial B_{r_1}}\rho^{2} (\varphi (x)) }
 \leq(\frac{ r_{2}}{ r_{1} })^{C},$$ when $0<r_{1}\leq r_{2}\leq  \varepsilon$.
 \begin{flushright}{\bf Q.E.D.}\end{flushright}
\section{The   adjunction numbers}
    For a real surface $\Sigma$ in 
a Riemannian 4-manifold $N$ which has an almost complex structure $J_N$, 
one can 
consider the intersection of $T_x\Sigma$ and
$J_{N}T_x\Sigma $ for points $ x \in \Sigma.$ 
There are only two possibilities: either $T_x\Sigma
\cap J_{N}T_x\Sigma =\{0\}$ where $x$ is called a totally real point or $T_x\Sigma
= J_{N}T_x\Sigma $ where $x$ is called a complex point. When the complex points are isolated,
it has a well-defined index at each complex point and there are formulas which relate  the 
total number of the complex points with indices to the topology of $\Sigma$. (See
\cite{CT},
\cite{ EH}, \cite{W1}, \cite{W2}, \cite{Wo}.)  The characterization given by J.~Chen and G.~Tian 
\cite{CT} is the following:
$$a_{N}(\Sigma)=\int _{\Sigma}(K_{T} +K_{N} ) \;dA=\sum ind\; x_{k},$$ 
where $K_T$ 
and $K_N$ are the Gaussian curvatures of the tangent bundle and  
normal bundle of $\Sigma$ in $N$ respectively and $ind\; x_k$ is the index at a 
complex point $x_k$. The first equality is the definition of  the adjunction 
number
 $a_{N}(\Sigma)$   of $\Sigma$ in $N$ and the second equality is a theorem
proved in \cite{CT}.
The tangent planes and normal planes on a branched minimal
 surface are still well defined even at branched points \cite{GOR}. The above
discussions
 also hold for branched minimal surfaces and in that case   the  
 integral is understood as
an improper integral. Moreover,
 it is proved by S.~Webster \cite{W1} and also by J.G.~Wolfson \cite{Wo}
that the complex points on a branched minimal surface are isolated 
and all of negative index when the surface is not holomorphic or 
antiholomorphic.

\VJ

\NO The bundle of   complex structures on $ R^{2l}$ along a minimal surface $\Sigma$
was discussed in  R.~Schoen's unpublished paper \cite{S2}. For the sake of completeness 
and the readers' reference, we adapt the argument to our settings and include a
discussion here.
One can identify the 2-vectors $\wedge ^{2}   R^{4}$ with the 
anti-symmetric $4 \times 4$ matrices by associating to a 2-vector $\eta$ 
$$\eta =\frac{1}{2} \sum  a_{kl}e_{k}\wedge e_{l}$$ 
the anti-symmetric matrix $A=(a_{kl})$, where $\{ e_{k},\; 1\! \leq \!k \!\leq \!4 \}$ is 
an oriented
orthonormal basis  of $  R^{4}.$  The inner product of $\wedge ^{2}   R^{4}$
induced on the anti-symmetric matrices is denoted by $<\cdot,\cdot>$, and it is 
$$<A,B>=-\frac{1}{2}Tr(AB)$$
for $A,B$ anti-symmetric matrices. Denote the set of oriented complex structures on 
$  R^{4}$ by $C_{4}$. That is, it is the set of positively oriented 
$ J:   R^{4} \rightarrow    R^{4}$ satisfying $$J^{t}J=I, \;\;\; J^{2}=-I,$$
where $J^{t}$ is the transpose of the matrix $J$. The image of $C_4$ under the above identification is
the sphere of radius $\sqrt{2} $ in $\,\wedge ^{2}_{+}  R^{4}\,$ which consists of the 
self-dual 2-vectors in $\wedge ^{2} R^{4}$. Let $\{ f_{k},\; 1\! \leq \!k \!\leq \!4 \}$  
be
another oriented
orthonormal basis  of $  R^{4},$ where $f_{k}=\sum  m_{lk}e_l$ and denote
$ M= (m_{kl}).$ Note
that $M^{t}M=I,$ and thus $M^{t}=M^{-1}.$ If a 2-vector
 $\eta$ is identified with a matrix $A$ in the basis  
$\{ e_{k},\; 1\! \leq \!k \!\leq \!4 \}$, it 
is identified with the matrix $M^{-1}AM$ in the basis $\{ f_{k},\; 1\! \leq \!k \!\leq \!4 \}$.   If   a complex structure in the basis $\{ e_{k},\; 1\! \leq \!k \!\leq \!4 \}$   is expressed as a matrix $J$, it is
 expressed as
   $M^{-1}JM$ in the basis  $\{ f_{k},\; 1\! \leq \!k \!\leq \!4 \}$. Thus  we   have the   identification as a bundle on a Riemannian 4-manifold $N$.
   Denote the total space of the restricted bundle on $\Sigma$ by ${\cal E}$.
 We claim that  ${\cal E}$ has an almost complex structure.
 The fiber $S^2$  has an almost complex structure    or we also can define the
 almost complex structure directly   from $C_4$ as follows. 
 Let $\cal A$ be the set of anti-symmetric $4 \times 4$ matrices.
For $J \in C_{4}$, one has $$T_{J}C_{4}=\{A\in {\cal  A} :AJ+JA=0 \}.$$ 
We  define the almost
complex structure    $$ {\cal  J} : T_{J}C_{4} \rightarrow T_{J}C_{4}$$
on $C_{4}$ by ${\cal J} (A)=AJ.\;$ It is easy to check that this is a right definition
 and it gives the almost complex structure  
on the fiber. The same construction
  gives  the almost complex structure on $C_{2l}$ for $l\!>\!2$ as well. 
Using the Levi-Civita connection we have a 
complement to the fiber which is called
a horizontal space and it can be 
identified with $T\Sigma$ via the projection map.  The identification induces 
an almost complex structure on the horizontal space. Therefore, we have the almost complex 
structure on the total space ${\cal E}$ and we will denote it still by $\cal J$. Assume that 
$u(t)$ 
is a section along a curve 
$\gamma(t)$ in $\Sigma$ and $\frac{d \gamma (t)}{d t} =T$. Then
 $(\gamma (t), u(t))$  is a
curve in ${\cal E}$ and the projection of the tangent vector
into   the fiber is just $\nabla_{T} u.$

\VJ

\NO Assume that $\{e_1, e_2, e_3, e_4\}$ is a local, oriented orthonormal 
basis of the tangent bundle $TN $ over $\Sigma$ such that $\{e_1, e_2\}$ is an oriented 
basis
of $T\Sigma$. We define an almost complex structure $J_\Sigma$ of $TN
$ along 
$\Sigma$ by
\begin{eqnarray*}
&J_\Sigma (e_1) = e_2, ~~~J_\Sigma (e_2)=-e_1,\\
&J_\Sigma (e_3) = e_4, ~~~J_\Sigma (e_4) = -e_3.
\end{eqnarray*}
 Hence $J_\Sigma$ is a section of the above bundle. J.~Chen and G.~Tian
 \cite{CT} shows that  
 $$ K_{T} +K_{N} =\Omega_{12}+\Omega_{34}+\frac{1}{2}|H|^{2}-\frac{1}{4}|\nabla
  J_{\Sigma}|^{2},$$ where 
  $\Omega_{kl} $ are some ambient curvatures, $H$ is the mean curvature on $\Sigma$ satisfying $$|H|^2=(h_{11}^{3}+h_{22}^{3})^{2}+ 
(h_{11}^{4}+h_{22}^{4})^{2}$$ and $$|\nabla
  J_{\Sigma}|^{2}=2(h_{12}^{4}-h_{11}^{3})^{2}+2(h_{12}^{3}+h_{11}^{4})^{2}
+2(h_{22}^{4}-h_{12}^{3})^{2}+2(h_{22}^{3}+h_{21}^{4})^{2}.$$ Thus one has
 $$a_{N}(\Sigma)=\int _{\Sigma}(\Omega_{12}+\Omega_{34}+\frac{1}{2}|H|^{2}-
 \frac{1}{4}|\nabla
  J_{\Sigma}|^{2})\;dA. $$
 We will consider   maps from $\Sigma$ into $N$  from now on. Hence
the surface on the above discussions 
 should be replaced by the image of a map. But   we  will use the
 same notation whenever there is no confusion.
 \begin{lemma} 
 Assume that $\; \varphi :\Sigma \rightarrow (N,g) $  is a branched minimal immersion.
 Then the map  $J_{\Sigma}: \Sigma \rightarrow
 {\cal E}$ is holomorphic.
 \end{lemma}
{\bf Proof:}

\medskip

\NO Assume that $x^{1}, x^{2}$ are the conformal coordinates near a point
$p$ on $\Sigma$ for the pull back
metric. Denote 
the complex structure by $j$ which satisfies  $$j\frac{\partial}{\partial x^{1}}=
\frac{\partial}{\partial x^{2}}\;\;\mbox{ and }\;\;j\frac{\partial}{\partial x^{2}}=-
\frac{\partial}{\partial x^{1}}.$$
Let $\{e_1, e_2, e_3, e_4\}$  be a local, oriented orthonormal 
basis of the tangent bundle $TN $ as described before in a neighborhood of 
$\varphi(p)$ on $\varphi(\Sigma)$. Therefore $$d \varphi(\frac{\partial}{\partial x^1})
=\mu e_{1}\;\;\mbox{ and }\;\;d \varphi(\frac{\partial}{\partial x^2})
=\mu e_{2},$$ where $\mu=|d \varphi(\frac{\partial}{\partial x^1})|=
|d \varphi(\frac{\partial}{\partial x^2})|.$
Note that $J_{\varphi(\Sigma)}$, which we will denote by $J_{ \Sigma }$ instead,
 can be identified with
$$-(e_{1}\wedge e_{2}+e_{3}\wedge e_{4}).$$  
 When $p$ is an unbranched point, we have
\begin{eqnarray*} 
&&\nabla_{e_1} (e_{1}\wedge e_{2}+e_{3}\wedge e_{4})\\ 
&=&\nabla_{e_1}  e_{1}\wedge e_{2}+e_{1}\wedge \nabla_{e_1}  e_{2}+
\nabla_{e_1}  e_{3}\wedge e_{4}+e_{3}\wedge \nabla_{e_1}  e_{4}\\ 
&=&(h_{11}^{3}e_{3}+h_{11}^{4}e_{4}) \wedge e_{2}
+e_{1} \wedge (h_{12}^{3}e_{3}+h_{12}^{4}e_{4})\\ 
&&+(- h_{11}^{3}e_{1}-h_{12}^{3}e_{2}) \wedge e_{4}
+e_{3} \wedge (-h_{11}^{4}e_{1}-h_{12}^{4}e_{2})\\ 
&=&(h_{12}^{3}+h_{11}^{4})e_{1}\wedge
e_{3}+(h_{12}^{4}-h_{11}^{3})e_{1}\wedge e_{4}
\\ 
&&+(-h_{11}^{3}+h_{12}^{4})e_{2}\wedge
e_{3}+(-h_{11}^{4}-h_{12}^{3})e_{2}\wedge e_{4} 
\\ 
&=& (h_{12}^{4}-h_{11}^{3})(e_{1}\wedge e_{4}+e_{2}\wedge e_{3})
+(h_{12}^{3}+h_{11}^{4})(e_{1}\wedge e_{3}-e_{2}\wedge e_{4}).
\end{eqnarray*}
A similar computation gives 
\begin{eqnarray*} &&
\nabla_{e_2} (e_{1}\wedge e_{2}+e_{3}\wedge e_{4})\\ 
&=&(h_{22}^{4}-h_{12}^{3})(e_{1}\wedge e_{4}+e_{2}\wedge e_{3})
+(h_{22}^{3}+h_{21}^{4})(e_{1}\wedge e_{3}-e_{2}\wedge e_{4}).
\end{eqnarray*}
The 2-vectors are identified with the anti-symmetric matrices, so  we
mix the notations sometimes. It can be checked that 
$$ e_{1}\wedge e_{4}+e_{2}\wedge e_{3} \in T_{J_\Sigma}C_{4}\;\;\mbox{  and  } 
\;\; e_{1}\wedge e_{3}-e_{2}\wedge e_{4} \in T_{J_\Sigma}C_{4}.$$
Furthermore,  
we have $${\cal J} (e_{1}\wedge e_{4}+e_{2}\wedge e_{3})=
e_{1}\wedge e_{3}-e_{2}\wedge e_{4}$$
\NO and 
$${\cal J} (e_{1}\wedge e_{3}-e_{2}\wedge e_{4})
= -(e_{1}\wedge e_{4}+e_{2}\wedge e_{3}).$$
 Because one has $h_{11}^{3}=-h_{22}^{3}$
and $h_{11}^{4}=-h_{22}^{4}$ on a minimal surface, it follows that 
\begin{eqnarray*}
& &{\cal J}\nabla_{e_1} (e_{1}\wedge e_{2}+e_{3}\wedge
e_{4})\\ 
&=&(h_{12}^{4}-h_{11}^{3})(e_{1}\wedge e_{3}-e_{2}\wedge e_{4})
+(-h_{12}^{3}-h_{11}^{4})(e_{1}\wedge e_{4}+e_{2}\wedge e_{3})\\ 
&=&(h_{12}^{4}+h_{22}^{3})(e_{1}\wedge e_{3}-e_{2}\wedge e_{4})
+( h_{22}^{4}-h_{12}^{3})(e_{1}\wedge e_{4}+e_{2}\wedge e_{3})\\ 
&=&\nabla_{e_2}
(e_{1}\wedge e_{2}+e_{3}\wedge e_{4}).
\end{eqnarray*} 
\NO That is, we have $\;{\cal J}\nabla_{e_1} J_{\Sigma}
=\nabla_{e_2} J_{\Sigma}\;$. Since the almost complex structure on the
horizontal space is given by the identification with $T\Sigma$, the map
also satisfies the holomorphic condition in the
horizontal space. Thus $$dJ_{\Sigma}  
(j \frac{\partial}{\partial x^1})=dJ_{\Sigma}  
(\frac{\partial}{\partial x^2})={\cal J}dJ_{\Sigma}  
(\frac{\partial}{\partial x^1}),$$ or  
$J_{\Sigma}$ is holomorphic away from the branched points. Because  $J_{\Sigma}$ is a
continuous map, it then follows that $J_{\Sigma}$ is in fact holomorphic
at all points
on $\Sigma$ by the standard fact in complex analysis or see the discussions below.
\begin{flushright}{\bf Q.E.D.}\end{flushright}
We would like to write the holomorphic condition in local coordinates and show that it is equivalent to satisfying a first order elliptic system.
Let $\{e_1, e_2, e_3, e_4\}$  be a local, oriented orthonormal 
basis of the tangent bundle $TN $. Then \begin {eqnarray*}
&E_{1}=e_{1}\wedge e_{2}+ e_{3}\wedge e_{4}&\\
&E_{2}=e_{1}\wedge e_{3}-e_{2}\wedge e_{4}&\\
&E_{3}=e_{1}\wedge e_{4}+ e_{2}\wedge e_{3}& 
 \end{eqnarray*} 
becomes a local basis for the bundle of self-dual 2-vectors. Assume that $\sum
u_{i}E_{i}$ is a section of this bundle on $\Sigma$.
The covariant derivative of the section in the direction $\frac{\partial}{\partial x^{1}}$ is   
\begin {eqnarray*}
\nabla_{\frac{\partial}{\partial x^{1}}}\sum u_{i}E_{i}=\sum \frac{\partial u_{i}}{\partial x^{1}}E_{i}
+\sum u_{j}<\nabla_{\frac{\partial}{\partial x^{1}}}E_{j},E_{i}>E_{i}.
 \end{eqnarray*} 
Note that $C_{4}$ is identified with a sphere of radius $\sqrt{2} $ in $\,\wedge ^{2}_{+}  R^{4}.$ Thus if the section lies in this subbundle,
one has $$\nabla_{\frac{\partial}{\partial x^{1}}}\sum u_{i}E_{i}=a_{1}\xi_{1}+
b_{1}\xi_{2} \hskip.3cm \mbox{and} \hskip.3cm 
\nabla_{\frac{\partial}{\partial x^{2}}}\sum u_{i}E_{i}=a_{2}\xi_{1}+
b_{2}\xi_{2},$$
where $\{\xi_{1},\xi_{2}\}$  satisfies
$${\cal J} (\xi_{1})=\xi_{2} \hskip.3cm \mbox{and} \hskip.3cm 
{\cal J} (\xi_{2})=-\xi_{1}.$$ A section is holomorphic is then equivalent to $a_{1}=b_{2} $ and $b_{1}=-a_{2}.$ Assume that   $J_{\Sigma}$ is identified with
$$-(e_{1}\wedge e_{2}+e_{3}\wedge e_{4})$$ at a point $p$ and is written as
$J_{\Sigma}=\sum u_{i}E_{i}$ near $p$, where $u_{1}=-\sqrt{1-u_{2}^{2}-u_{3}^{2}}.$
Then   \begin {eqnarray*}
a_{1}&=&<\nabla_{\frac{\partial}{\partial x^{1}}}\sum u_{i}E_{i}, \xi_{1}>\\
&=&\sum \frac{\partial u_{i}}{\partial x^{1}}<E_{i},\xi_{1}>+\sum u_{j}\Gamma_{1j}^{i}<E_{i},\xi _{1}>  \end{eqnarray*} 
and $b_{1}, \;  a_{2},\; b_{2}$  also have  similar expressions. We can choose 
$\xi _{1}=E_{2}$ and $\xi _{2}=-E_{3}$ at $p$. Then at $p$,
$$a_{1}=  \frac{\partial u_{2}}{\partial x^{1}}+\sum u_{j}\Gamma_{1j}^{2} \hskip1cm b_{1}= - \frac{\partial u_{3}}{\partial x^{1}}-\sum u_{j}\Gamma_{1j}^{3},$$ and 
$$a_{2}=  \frac{\partial u_{2}}{\partial x^{2}}+\sum u_{j}\Gamma_{2j}^{2} \hskip1cm b_{2}= - \frac{\partial u_{3}}{\partial x^{2}}-\sum u_{j}\Gamma_{2j}^{3}.$$ The symbols for the equations $a_{1}-b_{2}=0$
and $a_{2}+b_{1}=0$ are  nondegenerate at $p$. By continuity, they 
are still nondegenerate near $p$. Hence $u_{1}, u_{2}$ satisfy
a first order elliptic system and the coefficients of the lower order terms are
related to $\Gamma_{\alpha j}^{i}$ only. By an interior Schauder estimate \cite{ADN}, 
one has 
$$|J_{\Sigma}|_{1,\alpha ;B_{\varepsilon}} \leq C( |J_ {\Sigma}|_{0;B_{2 \varepsilon}}+|f|_{0,\alpha ;B_{2 \varepsilon} }),$$
where $C$ is a constant and $f $ is the zero order term.
 Hence if $|J_{\Sigma}|$ is bounded,   the isolated singularity is removable.
\section{The limit of the adjunction numbers}
\begin{theorem}
 Let $\varphi_{i} :\Sigma \rightarrow (N,g_i )$ be  a stable branched minimal 
 immersion    from a closed surface $\Sigma$ to a Riemannian  $4$-manifold 
$(N,g_i )$. Assume that $g_i$ converges to $g_0$ and $\varphi_{i}$ converges
 to 
   $\varphi_0$ in $C^{\infty}$, where $\varphi_0$ is a branched minimal 
immersion 
   from $\Sigma$ to $(N,g_0 )$.
 Then $$a(\varphi_{0} (\Sigma))=\lim_{i \rightarrow \infty} a(\varphi_{i} 
 (\Sigma)).$$
 \end{theorem}
\noindent  {\bf Proof:}

\medskip

\noindent  Without  loss of generality, we can assume that $\varphi_{0}$
has only one branched
 point at $x_0$. Let $B_r (x_{0})$ be the ball centered at $x_0$ 
 of radius $r$  with
 respect to the pull back metric $\varphi_{0} ^{\ast} (g_{0}).$ For $i$ large enough  
 all the branched points of $\varphi_{i}$ are within $B_r (x_{0})$ and $K_{T}^{i}+K_{N}^{i}$
 converges to $K_{T}^{0}+K_{N}^{0}$ on $\Sigma\setminus B_r (x_{0})$ uniformly. We have
 \begin{eqnarray*}
 a(\varphi_{0} (\Sigma))&=&\lim_{r \rightarrow 0} \int_{\Sigma\setminus B_r (x_{0})}
 (K_{T}^{0}+K_{N}^{0}) \;dA_{0}\\ 
 &=&\lim_{r \rightarrow 0} \lim_{i \rightarrow \infty} \int_{\Sigma\setminus 
 B_r (x_{0})} (K_{T}^{i}+K_{N}^{i}) \;dA_{i}\\ 
 &=&\lim_{i \rightarrow \infty} a(\varphi_{i} (\Sigma))-
 \lim_{r \rightarrow 0} \lim_{i \rightarrow \infty} \int_{B_r (x_{0})} 
 (K_{T}^{i}+K_{N}^{i}) \;dA_{i} ,
 \end{eqnarray*}
 where $dA_{i}$ is the volume form for the pull back metric $\varphi_{i} ^{\ast} (g_{i}).$  
 Because $$ K_{T}^{i} +K_{N}^{i} =\Omega_{12}^{i}+\Omega_{34}^{i}+
 \frac{1}{2}|H_{i}|^{2}-\frac{1}{4}|\nabla
  J_{i}|^{2}$$ and $H_{i}=0$ at unbranched points, we have
 \begin{eqnarray*}
 &&\lim_{r \rightarrow 0} \lim_{i \rightarrow \infty} \int_{B_r (x_{0})} 
 (K_{T}^{i}+K_{N}^{i}) \;dA_{i}\\ 
& =&\lim_{r \rightarrow 0} \lim_{i \rightarrow \infty} \int_{B_r (x_{0})}
 (\Omega_{12}^{i}+\Omega_{34}^{i}-\frac{1}{4}|\nabla J_{i}|^{2}) \;dA_{i}\\ 
 &=& \lim_{r \rightarrow 0}  \int_{B_r (x_{0})}
 (\Omega_{12}^{0}+\Omega_{34}^{0}) \;dA_{0}-
 \frac{1}{4}\lim_{r \rightarrow 0} \lim_{i \rightarrow \infty} 
 \int_{B_r (x_{0})}|\nabla J_{i}|^{2}  \;dA_{i}\\ 
& =&- \frac{1}{4}\lim_{r \rightarrow 0} \lim_{i \rightarrow \infty} 
 \int_{B_r (x_{0})}|\nabla J_{i}|^{2}  \;dA_{i}. 
 \end{eqnarray*}
 If we can show that $\lim\limits_{r \rightarrow 0} \lim\limits_{i \rightarrow \infty} 
 \int_{B_r (x_{0})}|\nabla J_{i}|^{2}  \;dA_{i}=0$, then the theorem will be proved. 
   Express the pull back metric as
 $h_{i}=\varphi_{i} ^{\ast} (g_i)=\lambda_{i} ^{2} \bar{h}_i,$  where 
$\bar{h}_i $ is a smooth metric with the volume form $d\bar{A}_i $
and $\lambda _i$ is a smooth scalar function 
with isolated and finite order zeros. We can choose $\lambda _i$ suitably such that $\lambda _i$ 
and  $\bar{h}_i $ converge to $\lambda _0$ 
and  $\bar{h}_0 $ in $C^{\infty}$ respectively.  Choose $r$ small enough such
 that $B_r (x_{0})$ is a conformal neighborhood for all $\bar{h}_i $. Compose 
$\varphi_{i}$ with a conformal transformation on  $B_r (x_{0})$ if necessary, 
  we can assume that $x^{1},  x^{2}$ are
  the conformal coordinates  for all $\bar{h}_i $. Because the image is minimal,
by the discussions in last section, it follows that 
$|\nabla _{ \frac{\partial}{\partial x^k}} J_{i}|$ is bounded for any fixed $i$.
 That is,  the energy density of $J_i$ with respect to the metric
$\bar{h}_i $
is bounded for any fixed $i$.
We change the metric on the domain to $\bar h _{i}$, but still use 
 the same notation $|\nabla J_{i}|$. If  $|\nabla J_{i}|$ is bounded
in $ B_r (x_{0})$ by a constant $c$ which is independent of $i$, then it follows
\begin{eqnarray*}
  \lim_{r \rightarrow 0} \lim_{i \rightarrow \infty} 
 \int_{B_r (x_{0})}|\nabla J_{i}|^{2}  d\bar A_{i} &\leq& c 
 \lim_{r \rightarrow 0} \lim_{i \rightarrow \infty}
 \int_{B_r (x_{0})}  d\bar A_{i}\\ 
&=&c \lim_{r \rightarrow 0} \int_{B_r (x_{0})}  d\bar A_{0}
 \\ 
&=&0.
 \end{eqnarray*} 
 When the domain is two dimensional, the energy is a conformal  
 invariant on the metric of the domain.
  Hence   the left hand side is exactly the quantity
 we want to control. So that in this case the theorem follows. 

\VJ
\NO  Now assume that $$\max_{x\in B_r (x_{0})}
 |\nabla J_{i}(x)|=b_{i} \hspace{3em} \mbox{ for } i>0,$$ where $b_i$ tends to $\infty$
 and assume that the maximun value $b_{i}$ is obtained at $x_i$. 
   Because
 $K_{T}^{i}+K_{N}^{i}$
 converges to $K_{T}^{0}+K_{N}^{0}$ uniformly on $\Sigma\setminus B_r (x_{0})$ for any $r$,
 the sequence $x_{i}$ must converge to $x_0$.
 We define a new metric $h_{i}'=b_{i}^{2} {\bar h}_{i}$ on the domain and choose a ball 
 of radius $\frac{b_{i}r}{2}$ around $x_i$ with respect to $h_{i}'$.   If we denote  the energy density with respect to $h_{i}'$ still by $ \nabla J_{i} $,
 then we have $|\nabla J_{i}(0)|=1$ and $|\nabla J_{i}(x)|\leq 1$ for $x \in 
 B_{\frac{b_{i}r}{2}}(0)$. 
 Because $J_i$   satisfies a first order
elliptic system in local coordinates, by an interior Schauder estimate \cite{ADN}, 
one has 
$$|J_{i}|_{1,\alpha ;B_1}\leq C( |J_{i}|_{0;B_2}+|f_{i}|_{0,\alpha ;B_2}),$$
where $C$ is a constant and $f_i$ is related to the Christoffel symbol of 
$h'_i$ only. (See the  discussions in the end of last section.) Since
the metric $h'_i$ converges to the flat metric on $B_2$, it follows that 
$|f_{i}|_{0,\alpha ;B_2}$ converges to~$0$. Thus $|J_{i}|_{1,\alpha ;B_1}$
is uniformly bounded. By the Ascoli-Arzela convergent Theorem, we have $J_i$
converges to a section $\bar{J}$ uniformly in $C^{1}$ and 
   \begin{eqnarray} |\nabla \bar{J} (0)|=\lim_{i \rightarrow \infty}
   |\nabla J_{i}(0)|=1.\end{eqnarray}
 
\VJ
\NO  Note that $B_{r}(x_{0})$ is a conformal neighborhood for $\bar{h}_{i}$
with conformal coordinates $x^{1}, x^{2}$. With the coordinates, we denote the
ball of radius $\frac{r}{2}$ at $x_{i}$ in the Euclidean metric by 
$D_{\frac{r}{2}}(0)$. The map $\varphi_{i}$ is a conformal harmonic map from
$(D_{\frac{r}{2}}(0),\sum_{\alpha =1}^{2} (dx^{\alpha})^{2})$ to $(N,g_{i})$.
Define $\tilde{\varphi}_{i}(x)=\varphi_{i}(\frac{x}{b_i})$. Then 
$\tilde{\varphi}_{i}( x)$ is a conformal harmonic map from 
 $(D_{\frac{b_{i}r}{2}}(0),\sum_{\alpha =1}^{2} (dx^{\alpha})^{2})$ to 
 $(N,g_{i})$. Let $\rho^{2} (y, g_{i})$   be the square of the distance   between
$y  $ and  $\varphi_{i} (0)$ in $(N,g_{i})$.
   Assume
  that 
 $$\max_{D_{\frac{1}{b_{i}}}  }\rho^{2} (\varphi_{i}(x), g_{i}) =
 \max_{D_{1}  }\rho^{2}(\tilde{\varphi}_{i} (x), g_{i})
 =c_{i}^{2}. $$
  Because $b_i$ tends to $\infty$ , it follows that
 $c_i$ tends to $0$  and $\rho_{i}^{2}
 (\tilde{\varphi}_{i} (x))$ is
  a subharmonic
 function on $D_{1}(0)$ for $i$ large enough.
 Thus the maximun value $c_{i}^{2}$ for $\rho_{i}^{2}(\tilde{\varphi}_{i} (x))$
 can be attained at $\bar{x}_{i} \in  
 \partial D_{1}(0)$.
 By choosing a new parametrization we can assume that $\bar{x}_i$ is fixed, say
at the point $q=(1,0)$. 
 The image 
 $\tilde{\varphi}_{i} (D_{1})$ is a branched minimal surface in $(N,g_{i})$.
   Because 
  $g_i$ converges to $g_0$,
   the monotonicity constant for branched minimal surfaces and the radius
   where the bound holds can be chosen uniformly. Therefore,  \cite{Si} $$\mbox{area}\:
   ( \tilde{\varphi}_{i} (D_{1}), g_{i})=\mbox{area}\:
   (  \varphi_{i} (D_{\frac{1}{b_i}}), g_{i})<c c_{i}^{2}.$$
Define a new metric $ g_{i}'=\delta ^{2}c_{i}^{-2} g_{i} $ on $N$, where $\delta$
 is a constant determined later. Let $||\nabla \tilde{\varphi}_{i}||^{2}$ be
 the norm of the energy density of $\tilde{\varphi}_{i}$
 with respect to the metric $g_{i}'$.  
  Then $$  \int_{D_1}||\nabla \tilde{\varphi}_{i}||^{2}\;dA 
 =2\;\mbox{area}\,( \tilde{\varphi}_{i} (D_{1}), g_{i}')<c \delta ^{2}.  $$    
  If we choose $\delta$ small enough, then there will be no energy 
  concentration and a subsequence, which is still 
denoted by $\tilde{\varphi}_i ,$   converges to  a smooth harmonic
map $\varphi $ from $(D_{1}(0),\sum_{\alpha =1}^{2} (dx^{\alpha})^{2}) $ 
to $(R^4, \sum_{k=1}^{4}(dy^{k})^{2})$ in $C^{\infty} $
  by a result of J.~Sacks and K.~Uhlenbeck \cite{SU}.
  Moreover, $$ \rho^{2}(\tilde{\varphi}  (q),  \sum_{k=1}^{4}(dy^{k})^{2}) )
  =\lim_{i \rightarrow \infty}
   \rho ^{2}(\tilde{\varphi}_{i} (q),  g_{i}')=\delta ^{2}.$$ 
Hence  $\tilde{\varphi} $ is a nonconstant map.
\VJ
\NO For any $L>1$ we claim that the energy $E(\tilde{\varphi}_{i}(D_{L}),g_{i}')$
is also unformly bounded. This follows from a modification of the proof
of Lemma~5.
\VJ
\NO {\bf A modification of Lemma~5:}
   Because 
  $g_i$ converges to $g_0$ and  $\varphi_{i}$ 
converges to
$\varphi_{0}$ in $C^{\infty}$, there exists a uniform $\varepsilon$ such that  
$\varphi_{i}(D_{\varepsilon})$ lies in a normal neighborhood of 
$\varphi_{i}(0)$ in $(N,g_i)$ and  $\rho^{2}(\varphi_{i} (x), g_{i})$ 
is subharmonic on 
$D_{\varepsilon}$. Moreover, we also have that $$|\nabla \varphi_{i} (x)|\hskip.3cm
\mbox{and} \hskip.3cm \frac{|\Gamma_{kl}^{m}(\varphi_{i} (x))|}{\rho_{i}(\varphi_{i} (x))}$$
are uniformly bounded on $D_{\varepsilon}(0)$. Thus the constant $\bar{c}_i$ in Lemma~5 can be 
chosen uniformly. Because $\rho(\varphi_{i} (x), g_{i})$ converges
to $\rho(\varphi_{0} (x), g_{0})$ and   $\max_{\partial D_{\varepsilon}}\rho(\varphi_{0} (x),g_{0})$
is positive,  it follows that $\max_{\partial 
D_{\varepsilon}}\rho(\varphi_{i} (x), g_{i})$ has a uniform positive lower
bound. The constant $C_i$ in Lemma~5 can then be chosen uniformly. In conclusion,
we show that there exist positive constants $\varepsilon$ and $C$ such that 
the maps $\varphi_{i}:( D_{\frac{r}{2}}(0), \sum_{\alpha=1}^{2} (dx^{\alpha})^{2}) \rightarrow (N,g_{i})$ 
satisfies  $$\max_{D_{r_2} }\rho^{2} (\varphi_{i}(x), g_{i})
   \leq (\frac{r_2}{r_{1}})^{C} \max_{D
 _{r_1} }\rho^{2} (\varphi_{i}(x), g_{i})
   $$ for any 
 $0<r_{1}\leq r_{2}\leq  \varepsilon$. Note that a constant conformal 
 factor on the metric of the target will not affect the conclusion. 
 Because 
 $ \frac{L}{b_i} \leq \varepsilon$ for $i$ sufficiently large, Lemma~5 can be applied 
 to $\tilde{\varphi}_{i}$ on $D_{L}$. Thus
  we show that
    $$\max_{D_L} \rho^{2}(\tilde{\varphi}_{i} (x), g_{i}')  \leq L^{C}
  \max_{D_1}\rho^{2}(\tilde{\varphi}_{i} (x), g_{i}') \leq L^{C}  \delta ^{2}.$$ Since $\tilde{\varphi}_{i}(D_{L})=\varphi_{i}(D_{\frac{L}{b_i}})$ for $i$ sufficiently large, the image 
  lies in  a ball in $(N,g_{i})$
 where the monotonicity
 formula holds.  The same argument as 
 above shows that $\mbox{area}\,(\tilde{ \varphi}_{i} (D_{L}), g_{i}')<c\,L^{C} \delta ^{2} $ \cite{Si}.
   That is, the energy $E(\tilde{\varphi}_{i}(D_{L}),g_{i}')
  <C_{L}$, where $C_{L}$ is a constant depending on $L$ only.
  Therefore, there exists a subsequence of $\tilde{\varphi}_i$, which is still 
denoted by $\tilde{\varphi}_i ,$  such that
$\tilde{\varphi}_i$ converges to a smooth harmonic map $\varphi $ from $D_{L}(0)$  to $(R^{4},\sum_{k=1}^{4}(dy^{k})^{2})$ except
 finite points \cite{SU}.
Choose a sequence $L_k$ which tends to $\infty$ and use the diagonal process to choose
a subsequence    which converges to  a smooth harmonic
map $\varphi $ in any compact set of $R^2$ except
 finite points \cite{SU}.
Here $\varphi$ is a harmonic map
from $ (R^{2},\sum_{\alpha=1}^{2} (dx^{\alpha})^{2})$  to $(R^{4},\sum_{k=1}^{4}(dy^{k})^{2})$. The bubbling phenomenon (\cite{PW}, \cite{SU}) does not affect our discussions, so we will not concern the issue here. Consider the variations of $\varphi$
which have compact supports and vanish near the branched points. Because there are no branched points of $\tilde{\varphi}_{i}$ in the support
of the variation  for $i$ large enough, the stablity of $\tilde{\varphi}_{i}$ implies
the stablity of $\varphi$ (for such variations).
  By the same reason 
as above, we can show that the area of $ \tilde{\varphi}_{i}(D_{L})$ is of   
quadratic growth by  
the monotonicity formula \cite{Si}. Thus  the area of  $\varphi ( R^{2})$
 is also of quadratic growth. It is a theorem of M.J.~Micallef  
that every complete  and of  quadratic area 
growth stable branched minimal surface   
in $R^4$  is holomorphic with respect
to some complex structure on $R^4$ (\cite{M}, \cite{M2}). (The stablity is for 
variations which have compact supports and vanish near the branched points.) In particular, it implies that $\nabla J =0,$ 
where $J$ is the section associated with $\varphi (R^{2})$.
The section $J_i$ converges to $J$ in $C^{\infty}$ on any compact set away
from the branched points. 
On the other hand, we   know that $J_i$ converges to $\bar{J}$ in 
$C^{1,\alpha}$ on $B_{1}$ by (2). Thus $J=\bar{J}$ on $B_1$
and
$$|\nabla J (0)|=\lim_{i \rightarrow \infty}|\nabla J_{i}(0)|=1.$$ It is 
a contradiction. Hence $|\nabla J_{i} (x)|$ is uniformly bounded with respect to ${\bar h}
_{i}$ and $g_i$. The theorem is then proved.
\begin{flushright}{\bf Q.E.D.}\end{flushright}

\section{The main theorem}
\begin{theorem} 
Assume that $(N,g_{0})$ is a Kahler-Einstein surface with the first Chern class   negative.
Let $[A] $ be a class in the second homology group $  H_{2}(N,Z)$, which
can be represented by a finite union of branched Lagrangian  minimal surfaces with 
respect to the metric $g_{0}$. Then with respect to any other metric in the connected 
component of
$g_{0}$ in the moduli space of Kahler-Einstein metrics, the class $[A]$ can also be 
represented 
by 
a finite union of branched Lagrangian minimal surfaces. 
\end{theorem}
{\bf Proof:}

\medskip

\NO Let $g$ be any metric in the connected 
component of
$g_{0}$ in the moduli space of Kahler-Einstein metrics. There
exists a smooth family of Kahler-Einstein metrics $g_t,\; 0\leq t\leq 1,$
satisfying $g_{1}=g$. A metric is said to have the property P if
the class $[A]$ can   be 
represented 
by 
a finite union of branched Lagrangian minimal surfaces with respect to this 
metric. 
Let $$T=\{\;t\;|\;t\in [0,1] \mbox{ and } g_t \mbox{
has the
property P}\}.$$ From the assumption of the theorem, we know that $T$ contains
$0$.
Now assume that $t_0$ belongs to $T$. That is,  the class  can be 
written as  $[A]=\cup_{1}^{n}[\varphi_{i}(\Sigma_{i})]$, where
$\varphi_{i}:\Sigma _{i} \rightarrow (N, g_{t_0})$ is a branched minimal 
immersion and the image is Lagrangian. We will deform each $\varphi_{i}$
separately. So now we only work on a single map 
$\varphi_{t_0}:\Sigma  \rightarrow (N, g_{t_0})$.
 It is strictly stable   by {Proposition 1}. 
 Thus by {Theorem 2} there exists a strictly
stable branched minimal
immersion   $\varphi _t$ from $\Sigma$   to $(N,g_t)$
for $|t-t_{0}| <\varepsilon$  and $\varphi_t$ converges to 
  $\varphi_{t_0}$ in $C^{\infty}$.
The Lagrangian surface $\varphi_{t_0}(\Sigma)$  satisfies 
 $a(\varphi_{t_0} (\Sigma))=0$.  
  Because the adjunction number is  an
integer and $$\lim_{t \rightarrow t_0} a(\varphi_{t} 
 (\Sigma))=a(\varphi_{t_0} (\Sigma))$$ by Theorem~3, it
 follows that  $a(\varphi_{t} (\Sigma))=0$.
Since the complex points on a branched minimal surface which is not holomorphic or antiholomorphic
  are isolated 
and  of negative index  
(see \cite{W1},  \cite{Wo}), it follows that $\varphi_{t} (\Sigma)$ is totally 
real. A totally real, branched minimal surface in a Kahler-Einstein surface 
with 
 $C_{1}<0$  is Lagrangian (\cite{CT}, \cite{Wo}). Thus $ \varphi_{t} 
 (\Sigma)$ is a branched Lagrangian  
 minimal surface. Because there are only finite maps, we can choose 
 $\varepsilon$
 such that each $\varphi _{i}$ has a deformation in $|t-t_{0}| <\varepsilon$.
  Hence  the class $[A]$ can be 
represented by a finite union of branched Lagrangian  minimal surfaces with 
respect to the metric $g_{t }$ for $|t-t_{0}| <\varepsilon$.  
That is, the set $T$ is open.
 
\VJ

\NO Consider a smooth family of branched minimal immersions 
$\varphi_{t}:\Sigma  \rightarrow (N, g_{t })$, $t_{0}\leq t<b$, which can be
thought as the maps obtained from the above local deformation. Denote the area of $\varphi_{t}(\Sigma)$
in $(N, g_{t})$ by $A(\varphi_{t}, g_{t})$ and $h(t,x)=
 \varphi_{t} ^{\ast} (g_{t})(x)$ with volume form $dA_{t}$.
 Because $\varphi_{t} $ is a branched minimal immersion, the pull back 
 metric $h(t, x)=\lambda (t,x) ^{2} \bar{h}_{t}(t,x)$  for some 
 smooth metric
$\bar{h}_t $   with the volume form $d\bar{A}_t $.
Then \begin{eqnarray*}\frac { d A(\varphi_{t},g_{t})}{dt}&=&\int_{\Sigma}
\sum_{i,j=1}^{2}
 h^{ij}(t,x)\dot{h}_{ij}(t,x) dA_{t}\\&=&\int_{\Sigma}
\sum_{i,j=1}^{2}
\bar{ h}^{ij}(t,x)\dot{h}_{ij}(t,x) d\bar{A}_{t}. \end{eqnarray*} 
Note that $$ h_{ij} (t,x)=\sum_{k,l=1}^{4}
  {g}_{kl} (t, \varphi_{t}(x))
\frac{\partial \varphi_{t}^{k}(x)  }{\partial x^i}
\frac{\partial \varphi_{t}^{l} (x) }{\partial x^{j}}.$$ 
By Leibniz's rule, one knows that 
$\dot{h}_{ij}(t,x)$ comes from two parts: one is fixing $g_{t}$ and varying
$\varphi_{t}$ and another  is fixing $\varphi_{t}$ and varying $g_{t}$.
Because $\varphi_{t} $ is a branched minimal immersion, the 
contribution of the terms which are obtained from fixing $g_{t}$ and varying
$\varphi_{t}$  
is zero. Thus we only need to consider the situation where $\varphi_{t}$ is fixed and only $g_{t}$ is varied.
  In this case $$ \dot{h}_{ij} (t,x)
=\sum_{k,l=1}^{4}
 \frac{\partial {g}_{kl} (s, \varphi_{t}(x))}{\partial s}|_{s=t}
\frac{\partial \varphi_{t}^{k}(x)  }{\partial x^i}
\frac{\partial \varphi_{t}^{l} (x) }{\partial x^{j}}.$$ 
For fixed $t$ and $x$, we choose the conformal coordinates for $\bar{h}_{t}$
at $x$ such that
$h_{ij}(t,x)=\delta_{ij}\lambda ^{2}$ or $\bar{h}_{ij}(t,x)=\delta_{ij} $.
We also choose the normal coordinates 
for $g_t$ at $\varphi_{t}(x)$ such that 
$g_{kl}(t, \varphi_{t}(x))=\delta_{kl}$.
Then at $(t,x)$  $$ h_{ij}(t,x)=\sum_{k=1}^{4}
\frac{\partial \varphi_{t}^{k}(x)  }{\partial x^i}
\frac{\partial \varphi_{t}^{k} (x) }{\partial x^{j}}=\delta_{ij}\lambda ^{2}.$$
Because $g_t$ is a fixed smooth family, the quantity
$|\frac{\partial {g}_{kl} (s, \varphi_{t}(x))}{\partial s}|$ has a uniform
bound $c$. Hence
\begin{eqnarray*}
\frac { d A(\varphi_{t},g_{t})}{dt} &=&\int_{\Sigma}
\sum 
\bar{ h}^{ij}(t,x)  
 \frac{\partial {g}_{kl} (s, \varphi_{t}(x))}{\partial s}|_{s=t}
 \frac{\partial \varphi_{t}^{k}(x)  }{\partial x^i}
\frac{\partial \varphi_{t} ^{l}(x) }{\partial x^{j}} d\bar{A}_{t}
  \\&\leq&
 c\int_{\Sigma}\sum_{k,l,i}|\frac{\partial \varphi_{t}^{k}(x)  }
 {\partial x^i}|
|\frac{\partial \varphi_{t}^{l} (x) }{\partial x^{i}}|  d\bar{A}_{t}
 \\&\leq&\frac{c}{2}\int_{\Sigma}\sum_{k,l,i}(|\frac{\partial \varphi_{t}^{k}(x)
   }{\partial x^i}|^{2}+
|\frac{\partial \varphi_{t}^{l} (x) }{\partial x^{i}}|^{2})  
d\bar{A}_{t}\\&=&
8\,c\int_{\Sigma} \lambda ^{2} d\bar{A}_{t}\\
&\leq& 8\,c\,
A(\varphi_{t},g_{t}).
\end{eqnarray*} 
Hence $A(\varphi_{t},g_{t})\leq  e^{8\,c(t-t_{0})}A(\varphi_{t_0},g_{t_0}).$ 
We get a  uniform bound for the area. The
Gauss equation for minimal surfaces is
 $$\overline{K}_{N}(t,x)=K_{\Sigma}(t,x)+|I\! I|_{t}^{2}(x)$$ at   unbranched points,
where  $\overline{K}_{N}(t,x)$ is the sectional curvature of $(N, g_{t})$ on the tangent plane of $ \varphi_{t}
(\Sigma)$ at $ \varphi_{t} (x)$, $K_{\Sigma}(t,x)$ is the Gaussian
curvature for the pull back 
 metric $h(t, x)$, and $|I\! I_{t}|^2(x)$ is the norm of the second fundamental form of 
$ \varphi_{t}(\Sigma)$ in $(N, g_{t})$ at $ \varphi_{t} (x)$.
Integrating both 
 sides of the equation,
we get 
\begin{eqnarray*} \int_{\Sigma}\overline{K}_{N}(t,x)\;dA_{t}&=&\int_{\Sigma}
K_{\Sigma}(t,x)
\;dA_{t} +\int_{\Sigma}|I\! I_{t}|^{2}(x)\;dA_{t}\\ 
&=&2\pi \chi(\Sigma)
+2\pi B(t)+\int_{\Sigma}|I\! I_{t}|^{2}(x)\;dA_{t},
\end{eqnarray*} where $B(t)$ is the total branched order of the map $\varphi_{t}$.
 Note that the  
 integrals in the formula are all understood as
 improper integrals. Because we have the area bound for $\varphi_{t}
(\Sigma)$  and  $\overline{K}_{N}(t,x)$ is bounded for a fixed family $g_{t}$, $0\leq 
t\leq 1,$ it follows
that $\int_{\Sigma}|I\! I_{t}|^{2}(x)\;dA_{t}+2\pi B(t)$ is uniformly
bounded.  
\VJ

\NO Similar to the lemma of H.I.~Choi and R.~Schoen in \cite{CS},  one can show
 the boundness of the sup norm of 
$|I\! I_{t}|^2(x)$ in a  ball away from the branched points, if   the $L^2$
norm of $|I\! I_{t}|^2(x)$ is sufficiently small in a bigger ball (see \cite{An}).
Since $\int_{\Sigma}|I\! I_{t}|^{2}(x)\;dA_{t}$ is uniformly bounded, by
applying J.~Sacks and K.~Uhlenbeck's covering argument
\cite{SU}, one can pick up a subsequence which converges to a
branched minimal surface. Because the surfaces are Lagrangian and the
areas   are bounded, the
limit surface is  Lagrangian (see \cite{SW} and compare with \cite{L2})
and the area of the limit surface is also
bounded by the same constant. The area of a closed minimal surface
has a lower bound which depends only on the injective radius of the ambinent manifold.
Because $g_{t}$ is a fixed smooth family of metrics for $0\leq 
t\leq 1,$ this lower bound can be chosen uniformly. Thus once we have the area
bound for the union of closed minimal surfaces, the total number of closed
minimal surfaces in that union will   be bounded. 
\VJ

\NO If $T$ contains $t_{0}$, then $T$ contains 
$(t_{0}-\varepsilon, t_{0}+\varepsilon)$  by the local deformation. The argument in last paragraph
shows that the end points also belong to $T$. Thus one can apply the local deformation to the end points and continue the process. 
  Although we do not
have a lower bound for the length of the interval where the local deformation holds. During the process,
we do have a global area bound for the  family of minimal surfaces obtained. Moreover, the topology
is bounded and the union is finite.
  The same argument as above shows the closedness of $T$. The nonempty
set $T$ is both open and closed. So it must be the whole set $[0,1]$. That is, the 
class $[A]$ can  be represented by 
a finite union of branched Lagrangian minimal surfaces with repect to the metric 
$g=g_1$.
This completes the proof.
\begin{flushright}{\bf Q.E.D.}\end{flushright}
Now we   give a simple application of the theorem.  Let $N=(M,g) \times (M,g),$ where $(M,g)$ is a closed Riemannian surface with a hyperbolic
metric. Assume that $f$ is a map from a closed surface $\Sigma$ to $M$ 
whose  induced map on the first foundamental group
 $\pi _1$ is injective. Then the  induced map  of $(f,f)$ on
 $\pi _1$   is also   injective and  
there exists a  branched minimal surface in the homotopy class of $(f,f)$
by a result of R.~Schoen and S.T.~Yau \cite{SY1}. 
Because the metrics on the two 
components of $N$ are the same, the branched minimal immersions in the 
homotopy class of $(f,f)$ must be of the form $(\bar{f},\bar{f})$
by the uniqueness of harmonic maps into a hyperbolic space \cite{Ht}.
Thus the branched minimal immersions are Lagrangian
 if we reverse the orientation on the second component. By the same argument as in \cite{Lee}, it follows that
the  branched minimal surface in the homotopy class is  unique since every branched minimal surface in the class
 is Lagrangian.
\VJ

\NO Now we change the metric on the second component in its
 moduli space of hyperbolic metrics. By Theorem~4 we still have 
the existence of the branched Lagrangian minimal surfaces in the homotopy class
with respect to the new  metric.
The Lagrangian minimal surfaces obtained here can be of different topology with 
the one we have 
in \cite{Lee}. One can also try to combine our existence result in \cite{Lee}  with {Theorem~4} to get the existence of the branched Lagrangian minimal 
surfaces in other classes.

\begin{flushright} 
N{\footnotesize{ATIONAL}} T{\footnotesize{AIWAN}} U{\footnotesize{NIVERSITY}}, T{\footnotesize{AIPEI}}, T{\footnotesize{AIWAN}}
\end{flushright}


\begin{thebibliography}{999}
\bibitem{ADN} S. Agmon \& A. Douglis \& L. Nirenberg, Estimates near the boundary for solutions of elliptic partial differential equations satisfying
general boundary conditions II, CPAM {\bf 17} (1964) 35-92.
\bibitem{An} M. Anderson, The compactification of a minimal submanifold in
Euclidean space by the Gauss map, preprint.
\bibitem{B} R.L. Bryant, Minimal Lagrangian submanifolds of  K\"{a}hler-Einstein
manifolds, Lec. Notes in Math. Vol. {\bf 1255} 1-12,  
Springer, Berlin Heidelberg New York, 
 1987.
\bibitem{C}  B.Y. Chen, Geometry of submanifolds and its application. Science 
University of Tokoyo, 1981.
\bibitem{CT}J. Chen \& G. Tian, Minimal surfaces in Riemannian 4-manifolds, 
 Geom. funct. anal. {\bf 7} (1997) 873-916.
\bibitem{CS} H.I. Choi \& R. Schoen, The space of minimal embeddings of a 
surface 
into a 3-dimensional manifold of positive Ricci curvature, Inv. Math. {\bf 81}
(1985) 387-394.
\bibitem{EH} Y. Eliashberg \& V. Harlamov, Some remarks on the number of 
complex points 
on a real surface in the complex one, Proceedings of Leningrad Int.
Topology Conference 143-148, 1982.
\bibitem{EL} J. Eells \& L. Lemaire,  Deformations of  metrics and associated 
harmonic
maps, Geometry and analysis  33-45, Indian Acad. Sci.,
Bangalore, 1980. 
\bibitem{EL2} J. Eells \& L. Lemaire, A report on harmonic maps, Bull. London
Math. Soc. {\bf 10} (1978) 1-68.
\bibitem{ES} J. Eells \& J.H. Samposon, Harmonic mappings of Riemannian manifolds, Amer. J. Math. {\bf 86} (1964) 109-160.
\bibitem{GOR} R. Gulliver \& R.Osserman \& H. Royden, A theory of branched immersions of surfaces, Amer. J. Math. {\bf 95} (1973) 750-812.
\bibitem{Ht} P. Hartman, On homotopic harmonic maps, Can. J. Math. {\bf 19}  (1967) 673-687.
\bibitem{HL} R. Harvey \& B. Lawson, Calibrated geometries, Acta Math. {\bf 
148} (1982) 
48-156.
 \bibitem{Law} B. Lawson, Lectures on minimal submanifolds. Vol. 1, Publish or 
Press, 
 Berkeley, 1980.
 \bibitem{L2} Y.I. Lee,   The Limit of Lagrangian Surfaces in $R^{4}$, Duke Math. J. {\bf 71}  (1993) 629-631.
 \bibitem{Lee} Y. I. Lee, Lagrangian minimal surfaces in K\"{a}hler-Einstein 
surfaces of 
 negative scalar curvature,  Comm. Anal. Geom. 
{\bf 2} (1994) 579-592 .
\bibitem{Mc} R.C. McLean, Deformations of calibrated submanifolds, preprint.
\bibitem{M} M.J. Micallef, Stable minimal surfaces in Euclidean space, J. Diff. 
Geom. {\bf 19} (1984)
57-84.
\bibitem{M2} M.J. Micallef, A note on branched stable two-dimensional minimal
surfaces, miniconference on geometry and partial differential
equations (Canberra, 1985), 157--162, Proc. Centre Math. Anal. Austral. Nat. Univ., 10, Austral. Nat. Univ., Canberra, 1986.   
\bibitem{MW} M.J. Micallef \& J. G. Wolfson, The second variation of area of
 minimal surfaces in four-manifold, Math. Ann. {\bf 295} (1993) 245-267. 
\bibitem{Mo} J.K. Moser, On the volume elements on manifolds, Trans. Amer. Math.
Soc. {\bf 120} (1965) 280-296.
\bibitem{O} Y.G. Oh, Second variation and stabilities of minimal Lagrangian 
submanifolds 
in K\"{a}hler manifolds, Inv. Math. {\bf 101} (1990) 501-519.
\bibitem{PW} T. Parker \& J.G. Wolfson, A compactness theorem for Gromov's
moduli space, J. Geom. Anal. {\bf 3} (1993) 63-98.
\bibitem{Si} L. Simon, Lectures on geometric measure theory. Proc. Centre Math. Anal. Austral. Nat. Univ., 3, Austral. Nat. Univ., 1983.  
\bibitem{SU} J. Sacks \& K. Uhlenbeck, The existence of minimal immersions of 
2-spheres, 
Ann. Math. {\bf 113} (1981) 1-24.
\bibitem{S2} R. Schoen,  Compactness, regularity, and almost holomorphicity 
results for
stable minimal surfaces in arbitrary codimension, unpublished. 
\bibitem{SW}R. Schoen \& J.G. Wolfson, Minimizing volume among Lagrangian
submanifolds, preprint.
\bibitem{SY1}R. Schoen \& S.T. Yau,   Existence of incompressible minimal 
surfaces and 
the topology of three dimensional manifolds with non-negative scalar curvature, 
Ann. Math. 
{\bf 110} (1979) 127-142.
\bibitem{SY2}R. Schoen \& S.T. Yau, Lectures on harmonic maps. Conference
Proceedings and Lecture Notes in Geometry and Topology Vol 2, International Press, 1997.
\bibitem{SYZ} A. Strominger \& S.T. Yau \& E. Zaslow, Mirror Symmetry is T-duality, Nuclear Phys. {\bf B 479} (1996) 243-259. 
 \bibitem{W1} S. Webster, Minimal surfaces in a K\"{a}hler surface, J. Diff. 
Geom.  {\bf 20} (1984) 463-470.
\bibitem{W2} S. Webster, On the relation between Chern and Pontrjagin numbers,
Contemporary Math. No. 49  135-143, Amer. Math. Soc., Providence, RI, 1986.
\bibitem{Wo} J.G. Wolfson, Minimal surfaces in K\"{a}hler surfaces and Ricci 
curvature, 
 J. Diff. Geom. {\bf 29} (1989) 281-294.
\bibitem{Wo2} J.G. Wolfson, Minimal Lagrangian diffeomorphisms and the Monge-Amp\'{e}re equation,  J. Diff. Geom.  {\bf 46} (1997) 335--373.
\end{thebibliography}
\end{document}